\documentclass[hyperref]{article}

\usepackage{xspace}
\usepackage{amsmath,amsfonts}
\usepackage[T1]{fontenc}
\usepackage{color}
\usepackage{hyperref}

\usepackage{graphicx}
\usepackage{ifpdf}
\ifpdf
\DeclareGraphicsExtensions{.pdf}
\fi

\newcommand{\code}[1]{\texttt{#1}}
\newcommand{\flag}[1]{\str{#1}\xspace}

\newcommand{\file}[1]{\code{#1}}
\newcommand{\newpackagename}[1]{%
	\expandafter\def\csname #1\endcsname{\code{#1}\xspace}%
}
\newpackagename{MATLAB}
\newpackagename{OCTAVE}
\newpackagename{BASH}
\newpackagename{EXPODE}
\newpackagename{EXPODEDEV}
\newpackagename{EXPINT}

\newcommand{\EXPODECMD}{\texttt{expode}\xspace}

\newcommand{\newmethod}[1]{%
	\expandafter\def\csname #1\endcsname{\code{\lowercase{#1}}\xspace}%
}
\newcommand{\EXP}[1]{{\texttt{exp#1}}\xspace}
\newmethod{EXPRK}%
\newmethod{EXPRB}%
\newmethod{EXPMS}%
\newmethod{EXPMSSEMI}%
\newmethod{EXPNEW}%
\newmethod{ODE}%
\newmethod{ODESET}%
\newmethod{EXPSET}%
\newmethod{EXPRBSET}%
\newmethod{EXPRKSET}%
\newmethod{OPTINFO}%
\newmethod{EXPRBINFO}%
\newmethod{DEVALEXP}%
\newmethod{DEVAL}%

\newcommand{\refOption}[1]{\code{#1}}%
\newcommand{\refOptionPrefix}[2]{\refOption{#2}}%
\newcommand{\str}[1]{\code{'#1'}}%
\newcommand{\optVal}[1]{\str{#1}}%
\newcommand{\C}{\mathbb{C}}
\newcommand{\R}{\mathbb{R}}

\newlength{\indentsize}
\def\setnewline#1{%
	\let\oldnewline\\%
	\renewcommand{\\}{\oldnewline[0cm]{}\normalfont\ttfamily#1{}\csname setnewline\endcsname{#1}}%
}
\def\resetnewline{%
	\let\\\oldnewline%
}
\newenvironment{sniplet}[1]{%
	\renewcommand{\c}[1]{\% ##1}%
	\newcommand{\var}[1]{##1}%
	\renewcommand{\str}[1]{'##1'}%
	\tabbing%
	\setnewline{#1}%
	\normalfont\ttfamily#1%
}{\resetnewline\endtabbing}

\newcommand{\ode}{ODE}

\newcommand{\uu}{\ensuremath{u}}
\newcommand{\UU}{\ensuremath{Y}}
\newcommand{\hh}{\ensuremath{h_n}}
\newcommand{\hhms}{\ensuremath{h}}
\renewcommand{\AA}{\ensuremath{J}}
\renewcommand{\gg}{\ensuremath{g}}
\newcommand{\DD}{\ensuremath{D}}
\newcommand{\GG}{\ensuremath{G}}
\newcommand{\FF}{\ensuremath{F}}
\newcommand{\JJn}{\ensuremath{J_n}}
\newcommand{\ggn}{\ensuremath{g_n}}
\newcommand{\ddn}{\ensuremath{d_n}}

\newcommand{\ee}{\ensuremath{\operatorname{e}}}
\newcommand{\dd}{\ensuremath{\operatorname{d}\!}}

\newcommand{\expodeDefi}[1]{{\itshape #1}}

\usepackage[left=3.5cm,right=3.5cm,top=2cm,bottom=2cm]{geometry}

\begin{document}

\def\ourtitle{\EXPODE{} -- Advanced Exponential Time Integration
  Toolbox for \MATLAB}

\markboth{Georg Jansing}{\ourtitle}

\title{\ourtitle}

\author{GEORG JANSING, Heinrich-Heine Universit\"at D\"usseldorf
}

\maketitle

\begin{abstract}
  We present a \MATLAB{} toolbox for five different
  classes of exponential integrators for solving (mildly) stiff
  ordinary differential equations or time-dependent partial
  differential equations. For the efficiency of such exponential
  integrators it is essential to approximate the products of
  the matrix functions arising in these integrators with vectors in a
  stable, reliable and efficient way. The toolbox contains options for
  computing the matrix functions directly by diagonalization or by
  Pad\'e approximation. For large scale problems, Krylov subspace
  methods are implemented as well.

  The main motivation for this toolbox was to provide general tools
  which on one hand allows one to easily run experiments with
  different exponential integrators and on the other hand to make it
  easily extensible by making it simple to include other methods or
  other options to compute the matrix functions. Thus we implemented
  the toolbox to be compatible with the \ode{} solvers included in \MATLAB.
  Most of the methods can be used with adaptive stepping.
\end{abstract}

%
%
%
%





\section{Introduction}
\label{sec:introduction}

In this paper we present our new \MATLAB{} toolbox \EXPODE{} for exponential integrators containing some of the most prominent integrators developed recently. Our code is based on the implementations provided in \cite{HocLS98} and \cite{HocOS09}.

We strongly recommend the reader to consult the recent and well-written review on exponential integrators by Hochbruck and Ostermann \cite{HocO10Acta}. Historical remarks can be found in \cite{MinW05}.

Exponential integrators are usually applied to stiff ordinary
differential equations where the stiffness is generated by a linear
term.
The non-linearity is usually assumed to be well approximated by polynomials. Typical applications are abstract or discretized parabolic or hyperbolic partial differential equations.
The unprecedented memory capacity of modern computers can handle correspondingly fine spatial discretizations.
Implicit time integration schemes need to solve non-linear systems of equations with a dimension which is a multiple of the dimension of the underlying differential equation.
As a consequence, the solution of these systems becomes very costly.
Classical explicit schemes, while generally less costly per time step, offer a different problem. They often have interdependent time step and spatial element
size requirements for stability. Explicit exponential methods can often avoid or at least weaken these
conditions to allow a more efficient use of the memory capacity.

Target equations of our toolbox are very high-dimensional systems where
management costs for smart memory management and arithmetics are justified
by the high computational costs.

We allow direct computation of actions of the arising matrix functions with
vectors or approximation of those with a Krylov method. We have adaptive
step size control for some of the methods and support both autonomous and
non-autonomous equations.

For an easy entry we are as compatible as possible with \MATLAB's \ode{}
toolbox. We also give fine grained control over the integration process
and have a focus on extensibility such that new functionality can be added
with only minimal changes to \EXPODE's code.

In some numerical experiments we benchmark \EXPODE's performance. First we use the well known van der Pol equation to test our adaptive step size choices and then use problems typical for studying exponential integrators \cite{HocO05SIAM,HocLS98}. In a more extensive experiment in the end we investigate \EXPODE's potential when applied to the hyperbolic Maxwell's equations with the spatial discretization provided by another specialized and well written \MATLAB package, the discontinuos Galerkin finite elements provided by \cite{hesthaven2008nodal}.

\section{Exponential Integrators}
\label{sec:exponential-integrators}

In this section we briefly describe the class of problems which can be
solved by using our toolbox and the methods which are implemented.

\subsection{Problem class}
\label{sec:problem-class}

We consider nonlinear initial value problems
\begin{equation}\label{eq:pre-linearized-ode}
	\uu'(t) = \FF\bigl(t,\uu(t)\bigr), \qquad \uu(t_0) = \uu_0
\end{equation}
where $\uu: \R \to \C^d$ and $\FF: \R \times \C^d \to
\C^d$. 

Exponential integrators are based on linearization of
$\FF$. One can usually distinguish two different types. The first one
uses a fixed linearization,
\begin{equation}\label{eq:semilin-ode}
  \uu'(t) = \AA \uu(t) + \gg\bigl(t, \uu(t)\bigr), \qquad \uu(t_0) = \uu_0,
\end{equation}
where $\AA\in \C^{d\times d}$ is an approximation to the Jacobian of $\FF$. Roughly speaking, an exponential integrator can be expected to be efficient, if $\AA$ contains the stiff part of $\FF$ and if $\gg$ is nice in the sense that $\gg(t,\uu(t))$ is a smooth function in $t$, where $\uu$ is the solution of \eqref{eq:pre-linearized-ode}.

The second option is to linearize in each time step,
\begin{subequations}
	\label{eq:linearized-ode}
	\begin{gather}
    \uu'(t) = \FF(t,\uu(t)) 
    = \JJn \uu(t) \, + \ddn t\, + \ggn(t,\uu(t)), \qquad \uu(t_n) = \uu_n, 
	\end{gather}
where
\begin{equation}
       \JJn  = \frac{\partial \FF}{\partial \uu}(t_n, \uu_n), 
       \quad \ddn = \frac{\partial \FF}{\partial t}(t_n,\uu_n),
       \quad \ggn(t,\uu(t)) = \FF(t, \uu(t)) - \JJn \uu(t) - \ddn t.
\end{equation}
\end{subequations}
For details on the construction and on the analysis of exponential integrators for time-dependent partial differential equations we refer to the review  \cite{HocO10Acta} and the references therein.

\subsection{One-step methods}
\label{sec:one-step-methods}

\subsubsection{Exponential Runge-Kutta methods}
\label{sec:runge-kutta-methods}


A general exponential Runge-Kutta scheme applied to \eqref{eq:semilin-ode} is of the form
%
\begin{subequations}\label{eq:RK-scheme0}
	\begin{align}
		\UU_{ni}
			& = \ee^{c_i \hh \AA} \uu_n 
      + \hh \sum_{j=1}^s a_{ij}(\hh \AA) \GG_{nj}, \qquad 1 \leq i \leq s \\
		\uu_{n+1}
	&  = \ee^{\hh \AA} \uu_n + \hh \sum_{i=1}^s b_i(\hh \AA) \GG_{ni}, 
	\end{align}
\end{subequations}
where 
\begin{equation}
  \label{eq:Gni-def}
      \GG_{ni} := \gg(t_n + c_i \hh, \UU_{ni}).  
\end{equation}
The coefficient functions $b_i (z)$ are linear
combinations of the \expodeDefi{entire functions}
\begin{equation}
	\label{eq:phi-functions}
	\varphi_k(z) = \int_0^1 \ee^{(1 - \tau) z} \frac{\tau^{k-1}}{(k-1)!} \dd\tau
\end{equation}
and $a_{ij} (z)$ of $\varphi_k(c_i z)$ respectively. The following recurrence formula is satisfied
\begin{equation}
	\label{eq:phi-recursion}
	\varphi_{k+1}(z) = \frac{\varphi_k(z) - \varphi_k(0)} z, \quad \varphi_0(z) = \ee^z.
\end{equation}
These functions play a very important role in exponential integrators,
and the efficient computation or approximation of $\varphi_k(\gamma
\hh \AA) \uu$ for $\gamma \in [0,1]$ is essential, see 
section~\ref{sec:matrix-functions} below.

If the following \expodeDefi{simplifying assumptions} on the
coefficient functions $a_{ij}$ and $b_i$ are satisfied
\begin{equation}
	\label{eq:simplifying-assumptions}
	\sum_{i=1}^s b_i(z) = \varphi_1(z), 
   \qquad \sum_{j=1}^s a_{ij}(z) = c_i \varphi_1(c_i z), \quad 1 \leq i \leq s,
\end{equation}
cf. \cite{HocO05SIAM},
the scheme \eqref{eq:RK-scheme0} can be reformulated as
\begin{subequations}\label{eq:Runge-Kutta-scheme}
	\begin{align}
       \UU_{ni} & = \uu_n + c_i \hh \varphi_1 (c_i\hh \AA) \FF(t_n,\uu_n)
                    + \hh \sum_{j=1}^s a_{ij} (\hh \AA) \DD_{nj},
                 \qquad 1 \leq i \leq s \\
      \uu_{n+1} & = \uu_n + \hh \varphi_1 (\hh \AA) \FF(t_n,\uu_n)
                    + \hh \sum_{i=1}^s b_{i} (\hh \AA) \DD_{ni},
	\end{align}
\end{subequations}
where
\begin{equation}\label{eq:dnj-non}
  \DD_{nj} = \gg(t_n + c_j h_n, \UU_{nj}) - \gg(t_n, \uu_n).
\end{equation}
Thus, all (internal) stages can be interpreted as corrected
exponential Euler steps.

In our package, we restrict ourselves to explicit schemes
($a_{ij}(z) = 0$ for $j \geq i$) satisfying
\eqref{eq:simplifying-assumptions}.
Note that this implies $D_{n1} = 0$, and the sum over the inner stages
in \eqref{eq:Runge-Kutta-scheme} actually starts with two.

The package includes implementations of the exponential Euler method
($s=1$), two two-stage schemes proposed by Strehmel and Weiner
\cite[Example 4.2.2]{StrW92} (cf.\ (2.39) and (2.40) in
\cite{HocO10Acta}) and two three-stage methods by Hochbruck and Ostermann
(cf.\ (5.8) and (5.9) in \cite{HocO05SIAM}), two methods with four
stages, namely the ETD4RK method by Cox and Mattews in \cite{CoxM02}
and the ETD4RK-B method by Krogstad in \cite{Kro05}
and a method with five stages (cf.\ (5.19) by Hochbruck and Ostermann
in \cite{HocO05SIAM}). 

For details on the convergence properties of these methods in a
framework of parabolic partial differentail equations we refer
to \cite{HocO05SIAM}.

Our implementation allows one to easily add other explicit
Runge-Kutta schemes, see section \ref{sec:running-expode-integrators}
below. We implemented these methods for constant step sizes only.

\subsubsection{Exponential Rosenbrock-type methods}
\label{sec:rosenbrock-type-methods}

If we apply an exponential Runge-Kutta scheme to the linearized equation \eqref{eq:linearized-ode} we obtain \expodeDefi{exponential Rosenbrock-type methods} \cite{HocOS09}. For an efficient implementation, these methods should also be reformulated such that most of the matrix functions are multiplied with vectors of small norm:
\begin{equation}\label{eq:dnj-non-ros}
  \DD_{nj} = \gg_n(t_n+c_j h_n,\UU_{nj}) - \gg_n(t_n,\uu_n),
\end{equation}
where
\begin{equation*}
 \gg_n(t,\uu) = \FF(t,\uu) - \JJn \uu - \ddn t.
\end{equation*}
This yields
\begin{subequations} \label{eq:rosei-reform-nonauto}
\begin{align}
  \UU_{ni} = \uu_n &+ h_n c_i \varphi_1(c_ih_n\JJn) \FF(t_n,\uu_n) \nonumber\\
  &  + h_n^2 c_i^2 \varphi_2(c_ih_n \JJn) \ddn 
  + h_n\sum_{j=2}^{i-1} a_{ij}(h_n\JJn) \DD_{nj},
      \label{eq:inner_stages_nonauto}\\
  \uu_{n+1} = \uu_n &+ h_n \varphi_1(h_n\JJn) \FF(t_n,\uu_n)
  + h_n^2 \varphi_2(h_n\JJn) \ddn
  + h_n\sum_{i=2}^s b_i(h_n\JJn)\DD_{ni}.\label{eq:final_stage_nonauto}
\end{align}
\end{subequations}
For autonomous problems, we have $\ddn=0$.

The toolbox contains the \expodeDefi{exponential Rosenbrock-Euler}
method, where $s=1$, and the methods {\tt exprb3} and {\tt exprb4}
that are proposed in section~5.1 in \cite{HocOS09}.  For all three methods
an error estimator is available, which makes it possible to use
variable step sizes. The latter two methods use an embedded scheme for
the estimation of the local error. For the exponential Rosenbrock-Euler
method we used the error estimator  described in \cite{CalO09}.

\subsubsection{{\tt \EXP4}}
\label{sec:exp4-integrator}


The \EXP4 scheme proposed in \cite{HocLS98} was the seed for many
activities on exponential integrators and matrix functions. With the
paper, the authors provided \MATLAB\ and C-codes of \EXP4, which can be used
easily. The \EXP4 scheme has two different error estimators and also
features a dense output formula to evaluate the numerical solution at
arbitrary times. Our \EXPODE package is actually inspired by the
original \EXP4 codes. In particular a lot of fine tuning for the
Krylov method for the approximation of the matrix functions (see
section \ref{sec:matrix-functions}) was motivated by this integrator.

\EXP4 can be interpreted as a special case of an exponential
Rosenbrock-type method, which uses the $\varphi_1$ function only.




\subsection{Multistep Methods}
\label{sec:multistep-methods}

\subsubsection{Exponential Adams methods}
\label{sec:exponential-Adams-methods}


Motivated by classical Adams methods \cite{HaiW96,Nor69} their exponential
counterpart, \expodeDefi{exponential Adams methods}, were constructed for
the solution of semilinear problems \eqref{eq:semilin-ode} in \cite{HocO10BIT}.
In contrast to classical methods, the interpolation is done for the
nonlinearity $g$ instead of the full right hand side $F$.

For a constant step size $h$, an exponential $k$-step Adams method 
has the form
\begin{align}
	\label{eq:exp-adams-scheme}
	\notag \uu_{n + 1}
	& = \ee^{-\hhms \AA} \uu_n + \hhms \sum_{j=0}^{k-1} \gamma_j(-\hhms \AA) \nabla^j \GG_n \\
	& = u_n + \hhms \varphi_1(-\hhms \AA) (G_n - \AA \uu_n) + \hhms \sum_{j=1}^{k-1} \gamma_j(-\hhms \AA) \nabla^j \GG_n
\end{align}
with coefficient functions
\begin{equation}
	\label{eq:multistep-coefficients}
	\gamma_j(z) = (-1)^j \int_0^1 \ee^{(1 - \tau) z} \binom{-\tau} j \dd \tau,
	\quad \text{where} \quad
	\binom {-\tau} j = \frac 1 {j!} \prod_{k = 0}^{j - 1} (-\tau - k).
\end{equation}
Here, 
\begin{equation*}
	\nabla^0 G_n = G_n, \qquad \nabla^{j+1} G_n = \nabla^j G_n - \nabla^j G_{n-1}, \quad j \geq 1
\end{equation*}
denote the \expodeDefi{backward differences}
for $G_j = g(t_j, \uu_j)$, $j = 1, ..., n$.  The coefficient functions
are linear combinations of the $\varphi$ functions
\eqref{eq:phi-functions}. An analysis of these methods is given in
\cite{HocO10Acta}. To start the $k$-step methods, 
we used the fixed point iteration proposed in \cite{HocO11}. 
Alternatively, we also provide the option to use an 
exponential Runge-Kutta method to compute the starting values.

The package contains exponential $k$-step Adams methods for $k=1, \ldots, 6$.


\subsubsection{Linearized exponential multistep methods}
\label{sec:exponential-linearized-multistep}

The same idea of interpolation is now applied to the linearized
equation \eqref{eq:linearized-ode}. An additional order of accuracy is
gained by exploiting 
\begin{equation*}
  \frac {\partial g_n}{\partial \uu}(t_n, \uu_n) = 0,
  \qquad 
  \frac {\partial g_n}{\partial t}(t_n, \uu_n) = 0.
\end{equation*}
See \cite{HocO11} for details.
The interpolation polynomial $\widehat p_n$ now additionally satisfies $\widehat p_n'(t_n) = 0$ accordingly.
These
\expodeDefi{linearized exponential multistep methods} are defined as
\begin{subequations}  \label{eq:exp-multi-lin}
\begin{align} \label{eq:multi-new}
\uu_{n+1}
& = \uu_n + \hhms \varphi_1(\hhms \JJn) \FF(t_n, \uu_n) + \hhms^2 \varphi_2(\hhms \JJn) \ddn
+ \hhms \sum_{j=1}^{k-1} \widehat\gamma_{j+1} (\hhms \JJn) \sum_{\ell=1}^j \frac{1}{\ell} \nabla^\ell \GG_{n,n},
\end{align}
with weights
\begin{equation}  \label{eq:multi-weight-lin}
\widehat\gamma_{j+1}(z)  = (-1)^{j+1} \int_0^1 \ee^{(1-\tau)z}
\tau \binom{-\tau}{j} \dd \tau
\end{equation}
\end{subequations}
and backwards differences now based on $\GG_{n,m} = \gg_n (t_m , \uu_m)$ keeping the first index fixed.
In addition to this general scheme, an implementation of the
scheme proposed
in \cite[formula (39)]{Tok06} is contained in our toolbox.
%
%
Again using the fixed point iteration proposed in \cite{HocO11}, we implemented integrators for  $k = 1, \ldots, 5$.
Using an exponential Rosenbrock scheme for the initial steps instead, we can obtain exponential linearized multistep methods up to $k = 4$ due to the accuracy of the onestep methods. 



%

\section{Implementation Issues}
\label{sec:implementation-issues}

In this section we discuss some details of our
implementation. 


\subsection{Step Size Control}
\label{sec:step-size-control}

Step size control is provided for exponential Rosenbrock-type methods
and \EXP4\ via a standard Gustafsson approach \cite[pp.~31--35 and
pp.~550ff]{HaiW96} together with different norms of the scaled error
vector
\begin{equation*}
	\widehat e_n = \left( \frac {e_n(i)}{sc(i)} \right)_{i=1}^d,
	\qquad sc = ATol + \max\left\lbrace|\uu_n|, |\uu_{n-1}|\right\rbrace \cdot RTol,
\end{equation*}
where $(i)$ denotes the $i$-th component and $e_n$ is the estimated error in the $n$th time step. By default we use the maximum-norm (complying to \MATLAB{}'s defaults), but can easily switch to other norms like the Euclidian norm or user defined norms defined in a \MATLAB{} function or by an inner product defined via its Gramian matrix. This is implemented in the options \refOption{NormControl} and \refOption{NormFunction}.

The desired accuracy is determined by the
absolute ($ATol$) and relative ($RTol$) tolerance and the chosen
norms. An implementation using an iterative process for the 
evaluation of the matrix functions also has an impact on  the step size
selection. We will discuss this in the following
section.

\subsection{Matrix Functions}
\label{sec:matrix-functions}

We now turn our attention to the matrix functions arising in
exponential integrators. 

Since the $\varphi$ functions are analytic, $\varphi(h\AA)$ can be
computed via diagonalization if $\AA$ is diagonalizable.  The evaluation
of $\varphi$ at the eigenvalues can be computed directly from
\eqref{eq:phi-functions} via
\begin{equation}
  \label{eq:phi-direct}
  \varphi_k(z) = 
  \begin{cases}
    \frac{\ee^z - \left( 1 + z + \frac{z^2}{2!} + \ldots +
      \frac{z^{k-1}}{(k-1)!} \right)}{z^k(k-1)!}, & z \neq 0,\\
    \frac{1}{k!}, & z = 0.
  \end{cases}
\end{equation}
In finite precision arithmetic, we use a sufficiently high-order
Pad\'e approximation in a neighborhood of zero. Alternatively, one
could break down recursion formula \eqref{eq:phi-recursion} to the
evaluation of the exponential function $\varphi_0(z) = \ee^z$ or use
contour integration \cite{KasT05,SchT07,Lop10,TreWS06}.

Diagonalization is inefficient for large matrices, but fortunately,
exponential integrators do not require the complete matrix function
but $\varphi_k(\hh \AA) v$ for a vector $v$ only. This can be
approximated within a Krylov subspace with respect to $\AA$ and $v$,
see \cite{GalS92,HocL97}. 

Krylov methods have the advantage that they require the evaluation of
the matrix-vector products $\AA v$ only. It is not necessary to compute
the full Matrix $\AA$ explicitly. The error is controlled via the 
error estimators proposed in \cite{HocLS98,Saa92}.

We implemented configurable maximal Krylov subspace dimensions (cf. KrylovMaxDim option). If the Krylov approximation fails to converge within the maximum allowed dimension of the Krylov space, the step size $h$ has to be reduced such that with the reduced step size, the error estimator fulfills the accuracy requirements. If we cannot reduce the step size (e.~g.\ if we run a code for constant step sizes), a warning is triggered. Note that all products of $\varphi$ functions with the same vector $v$ can be approximated in one Krylov subspace. The dimension of this subspace is chosen such that all $\varphi$ functions are approximated sufficiently well.

In addition, we reuse data from previous steps if possible. For instance, if we compute matrix functions via diagonalization, solving semilinearized \ode{}s \eqref{eq:semilin-ode} and use
the same step size in the next time step, we do not have to recompute $\varphi_k(\hh \AA)$.
We also reuse the Krylov subspace for $\varphi_k(\hh \JJn) F(t_n, \uu_n)$ if we have to reduce the step size to to a step rejection.

For statistical purposes, we provide data about the Krylov process.
This data is displayed at the end of the integration process by default.
A typical output looks like this:
\begin{sniplet}{}%
	statistics: \\%
	[ ...~] \\%
	number of matrix function evaluation times vector:~331 \\%
	number of Krylov subspaces:~~~~~~~~~~~~~~~~~~~~~~~~147 \\%
	total number of Krylov steps:~~~~~~~~~~~~~~~~~~~~~~1455 \\%
	number of step size reductions due to Krylov:~~~~~~0 \\%
	number of recycled subspaces:~~~~~~~~~~~~~~~~~~~~~~10 \\%
	maximal dimensions of subspaces:~~~~F1:~15, v:~15,~D2:~15, D3:~11 \\%
\end{sniplet}

The first line reports the number of products of a $\varphi$-function
with a vector, the second line contains the number of Krylov
subspaces that have been built for different vectors. The third line counts the sum of all
Krylov subspace dimensions in the whole integration process. Next we
have the number of step size reductions which are necessary to fulfill
the error tolerances, followed by the number of reused
subspaces. Those can arise after a step rejection -- either by
the error estimator of the integrator or by the Krylov process. The
maximal dimensions in the last line correspond to the maximal
dimension of a subspace for a specific vector. Their meanings can be
looked up in table \ref{tab:KrylovStatistics}. Note that the multistep
methods require initial steps, such that -- depending on the choice of
their computation -- some of the labels for Runge-Kutta and
Rosenbrock-type methods may appear in their statistics as well.
\begin{table}[tbh]
	\begin{center}
	\begin{tabular}{l|l|l}
		name
			& vector
			& integrator type(s) \\
		\hline\hline
		\code{F1}
			& $\FF(t_n, \uu_n)$
			& linearized \\
		\code{v}
			& $\frac{\partial \FF}{\partial t}(t_n, \uu_n)$
			& linearized \\
		\code{D$\iota$}, $\iota = \texttt 2, \texttt 3$
			& $\DD_{n\iota}$
			& Rosenbrock-type \\
		\code{Y$\iota$plusA}, $\iota = \texttt 1, \dots, \texttt 5$
			& $\GG_{n\iota} - \AA \uu_n$
			& Runge-Kutta \\
		\code{d$\iota$}, $\iota = \texttt 4, \texttt 7$
			& $d_\iota$
			& \EXP4 \\
		\code{GDiff$\iota$}, $\iota = \texttt 1, \dots, \texttt 6$
			& $\nabla^\iota \GG_n$
			& both multistep \\
		\code{GDiffInit$\iota$}, $\iota = \texttt 1, \dots, \texttt 6$
			& $\Delta^\iota \GG_n$
			& both multistep\\[1ex]
	\end{tabular}
	\caption{\hspace*{-3cm} \text{Meaning of the labels in the Krylov statistics output}}
	\label{tab:KrylovStatistics}
	\end{center}
\end{table}




\section{Usage}

In this section we give a brief introduction to the \MATLAB package
\EXPODE. A more extensive documentation is contained in the \EXPODE
package (\code{manual.pdf}).



\subsection{Installation and Requirements}

We will now describe the minimal requirements and the installation of
the \EXPODE toolbox. \EXPODE runs on all recent and middle-aged
computers. The performance strongly depends on the problem and on the available
hardware. The toolbox was tested on \MATLAB versions down to \MATLAB
7.2 (R2006a), released in 2006. 
It is not compatible to versions prior to 7.0 due to the lack of
proper function handles. 
You can download the package from the author's web pages.
Two different versions are available, a package for users and an extended
one for developers. The latter contains some additional tools helpful
for extending \EXPODE. Usually the user package should be sufficient.

To install, just unpack the archive.
This will create a new \path{expode} subdirectory. To make it available
in \MATLAB, just add the package's root to \MATLAB's path and run the
\code{initPaths} function with
\begin{sniplet}{>{>}}
	addpath \path{mydownloadpath/expode}; \\
	initPaths;
\end{sniplet}
To make a permanent installation for the current user, put the above
line into your \file{startup.m} file. See \MATLAB's help for more
information.
%

\subsection{Quick Start}
\label{sec:quickstart}

To get a first impression of \EXPODE we start with running some of the
examples included. 
To access the examples we add the examples directory to the \MATLAB
path. Run the following commands
\begin{sniplet}{>{>}}
	addpath \path{mydownloadpath/expode/examples}; \\
	[\var{t}, \var{y}] = Heat1D([], [], 'run');
\end{sniplet}
to solve a heat equation with a time-dependent source term in one dimension. Use
\begin{sniplet}{>{>}}
	help Heat1D;
\end{sniplet}
to obtain information on the example. The solution will be visualized
in a mesh plot. All examples contained in the package
can be run by simply calling them without arguments. Short information
on the problems is available via \code{help}. To work with
the solver, it is convenient to run the example manually.
\begin{sniplet}{>{>}}
	\c{parameters} \\
	\var{N} = 100; \var{epsilon} = 0.1; \var{gamma} = 0.1; \\
	\\
	\c{get initial conditions} \\
	[\var{tspan}, \var{y0}, \var{options}] = Heat1D([], [], \str{init}, \var{N}, \var{epsilon}, \var{gamma}); \\
	\\
	\c{run the example} \\
	[\var{t}, \var{y}] = \EXPODECMD{}(@Heat1D, \var{tspan}, \var{y0}, \var{options});
\end{sniplet}
Now we can start playing with options and parameters. Switching to the
direct solver for the matrix functions, we use
\begin{sniplet}{>{>}}
	\var{options} = \EXPSET{}(\var{options}, \str{MatrixFunctions}, \str{direct});
\end{sniplet}
Other options are set similarly. 
Some checks on the values set for an option are applied automatically. 
An overview of the available options for an integrator is provided by
calling the integrator info without arguments. More detailed
information on a specific option can be shown with this command as
well:
\begin{sniplet}{>{>}}
	\EXPRBINFO~~~~~~~~~\c{prints all available options for exprb} \\
	\EXPRBINFO MinStep \c{prints helptext for MinStep option}
\end{sniplet}

A common task is to create order plots, where the problem is solved on
a fixed time interval with a number of different time steps and the
error is plotted over these time steps. Computing the error or an
approximation to it requires one to evaluate the exact solution or a very
accurate reference solution first. If an exact solution is available,
this can be done by calling \code{ode(t, [], 'exact')}. As an example
to show how simply this can be done with the package, we consider the
\code{semi1} example. We refer to its helptext for detailed information.
An order plot for a finite difference spatial discretization
with $\code{N} = 50$ grid points for all Rosenbrock-type methods is
created via
\begin{sniplet}{>{>}}
	allMethods(@semi1, \str{exprb}, \str{{}}, [], 50);
\end{sniplet}

The input argument \code{''} chooses the direct solver for the
evaluation of the matrix functions. The chosen step sizes depend on the
problem's \texttt{tspan} data and are chosen uniformly logarithmically.
This choice and other parameters can be manipulated with options to
\texttt{allMethods}.

To implement a new differential equation, we recommend modifying
one of the example files in the \path{examples} directory. You should
start in \path{examples/Hello_World}, where we put some introductory
files. \code{MinEx.m} is a very simple example while \code{Template.m}
uses more advanced features. Both files contain many helpful comments.

\subsection{Running Specific Integrators}
\label{sec:running-expode-integrators}
Here we briefly describe how the specific exponential integrators can be invoked. The calling sequence of all \EXPODE integrators is 
	\begin{sniplet}{>{>}}
		[t, y] 
        = integrator(@ode, \{@jac\}, \{tspan, y${}_0$\}, \{opts\}, \{varargin\});
	\end{sniplet}
	where \code{integrator} has to be substituted by either \code{expode}
        or one of the specific integrators below. Arguments in braces
        (\code{\{\}}) are optional. The at sign (\code{@}) represents
        either a function handle, a function name as string or an
        \code{inline} object. The function \code{@ode} has to evaluate
        data of the differential equation required for the solution. It
        follows \MATLAB's standard syntax, though it should be callable
        with
	\begin{sniplet}{>{>}}
		res = ode(t, y, flag, \{varargin\});
	\end{sniplet}
	where the flag controls what the function returns.
	The \code{ode} function has to return the evaluation of the right hand side
	of the differential equation, when the empty string (\flag{{}}) is given as
	flag. This is needed for all solvers. In addition other flags have
	to be handled depending on the integrator, see table \ref{tab:ODEFlags}.
	
\begin{table}[tbh]
	\begin{center}
	\begin{tabular}{p{0.15\textwidth}|p{0.27\textwidth}|p{0.27\textwidth}}
		flag
			& meaning
			& integrator type(s) \\
		\hline\hline
		\code{''}
			& $\FF(t, \uu)$
			& all \\
		\code{jacobian}
			& $\JJn = \frac{\partial \FF}{\partial \uu}(t, \uu)$
			& linearized \\
		\code{linop}
			& $\AA$ from \eqref{eq:semilin-ode}
			& semilinear (req.) and \newline
			  linearized (opt.) \\
		\code{gfun}
			& $\gg(t, \uu)$ from \eqref{eq:semilin-ode}
			& semilinear (req.) and \newline
			  linearized (opt.) \\
		\code{df\_dt}
			& $\ddn = \frac{\partial \FF}{\partial t}(t, \uu)$
			& linearized (opt.) \\
		\code{dg\_dy}
			& $\frac{\partial \gg}{\partial \uu}(t, \uu)$
			& semilinear (req.) and \newline
			  linearized (opt.) \\
		\code{init}
			& return default \newline
			\code{[ tspan, y0, opts ]}\newline
			for the equation
			& all (opt.) \\
	\end{tabular}
	\caption{\hspace*{-3cm} \text{Flags for the \code{ode}-file for the different integrators}}
	\label{tab:ODEFlags}
	\end{center}
\end{table}
	
	
	The \code{@jac} argument is only available for the linearized
        integrators and is a handle to a function evaluating the
        Jacobian. Alternatively the \code{ode} will be queried with
        the \flag{jacobian} flag. \code{tspan = [ t${}_0$, T ]} is the
        integration interval, \code{y${}_0$} the initial
        condition. \code{opts} is an options structure, set with one
        of the \code{set} commands and \code{varagin} will be passed
        to \code{ode} to configure parameters of the differential
        equation.
	
	The command to use an \emph{exponential Runge-Kutta
          integrator} is \EXPRK. To select one of the schemes
        described in section \ref{sec:one-step-methods} use the
        \refOptionPrefix{exprk}{Scheme} option. The default scheme is
        \optVal{Krogstad}. To set the parameters appearing in some of
        the methods, use the \refOptionPrefix{exprb}{Parameters}
        option. For a detailed overview of the available schemes we
        refer to the integrator documentation.
	
	The command to use an \emph{exponential Rosenbrock-type
          integrator} is \EXPRB. To select one of the three available
        schemes, use the \refOptionPrefix{exprb}{Order} option with
        value \optVal{two}, \optVal{three} or \optVal{four}, where the
        latter is the default.
	To control the parameters for the error estimator for the
        order four integrator use
        \refOptionPrefix{exprb}{ErrorEstimate}. We refer to the
        documentation for details on the schemes.
	
	The \EXP4\ method is called with the \code{exp4} command. It
        has a built in dense output generator, that allows one to
        evaluate the numerical solutions at arbitrary times. Specify
        \code{t = [ t${}_0$, \dots, t${}_m$ ]} to use this feature,
        where the sequence $\{t_j\}_{j=0}^m$ is either strictly
        increasing or strictly decreasing.
	
	The \emph{semilinear $k$-step methods} are available via
        \code{expmssemi}, where $k$ 
        is set with the
        \refOptionPrefix{expmssemi}{kStep} option. Our implementation
        allows constant step sizes only.
	
	The \emph{linearized $k$-step methods} are called with
        \code{expms}, where $k$ is defined as for \code{expmssemi}.
        You can also select \code{'Tokman'} to  use the scheme suggested in \cite[eq.~(39)]{Tok06}. 
	

\subsection{Properties of Equations}
Some problems allow one to exploit certain properties to improve the
efficiency of the integrator. In table \ref{tab:ode-properties} we
collected some properties together with information on how to exploit
them in \EXPODE. Note that for non autonomous equations you might get
wrong results, when the appropriate option is not set.
\begin{table}
\begin{center}
\renewcommand{\arraystretch}{1.2}
	\begin{tabular}{p{0.34\textwidth}|l|p{0.21\textwidth}|l}
		\ode property
			& option
			& choices
			& integrators \\
		\hline\hline
            autonomous/nonautonomous
			& \refOption{NonAutonomous}
			& \texttt{\{'off'\}}, \texttt{'on'}
			& linearized \\
                semilinear equation \eqref{eq:semilin-ode}
			& \refOption{Semilin}
			& \texttt{\{'off'\}}, \texttt{'on'}
			& linearized \\
		constant/nonconstant \newline
		Jacobian
			& \refOption{JConstant}
			& \texttt{\{'off'\}}, \texttt{'on'}
			& linearized \\
		complex/real valued solution
			& \refOption{Complex}
			& \texttt{'off'}, \texttt{\{'on'\}}
			& all \\
		structure of the Jacobian \newline
		or linear part
			& \refOption{Structure}
			& \texttt{\{'none'\}}, \newline
			  \texttt{'normal'}, \newline
			  \texttt{'symmetric'}, \newline
			  \texttt{'skewsymmetric'}, \newline
			  \texttt{'diagonal'}
			& all \\[1ex]
	\end{tabular}
	\caption{\hspace*{-3cm}\text{Options corresponding to properties of differential equations.}\newline
	\hspace*{-3cm}\text{Default values are set in braces}}
	\label{tab:ode-properties}
\end{center}
\end{table}

\subsection{Matrix Functions}

As mentioned before, the evaluation of the matrix functions is a crucial point in the implementation of exponential integrators.


The default setting for the matrix function evaluator is to compute
and store the full matrix functions. This is suitable for small or
medium sized problems.  It requires one to evaluate the full Jacobian or
its linear part. For linearized problems, the Jacobian has to be
computed in each time step while for methods based on a fixed
linearization, it has to be computed only once.

We also provide an implementation of the 
Arnoldi process  to approximate the matrix
functions, which can be used for large scale problems. 
It should not be used for very small problems, where the direct
evaluation is more efficient and more reliable.

To switch on the Krylov method, set 
\begin{sniplet}{>{>}}
	options = \EXPRBSET{}('MatrixFunctions', 'arnoldi');
\end{sniplet}
Krylov subspace methods can be implemented by using the matrices $\JJn$
or $\AA$ explicitly (saved as sparse matrices) or in a matrix-free
fashion, where subroutines for the evaluation of the matrix-vector
products $\JJn v$ or $\AA v$, respectively, are provided.  The options for
these matrix-free versions are activated by setting
\begin{sniplet}{>{>}}
	options = \EXPRBSET{}(options, 'JacobianV', 'on'); \\
	options = \EXPRBSET{}(options, 'LinOpV', 'on');
\end{sniplet}
respectively. If one of these options is \code{'on'}, then the
corresponding \ode file should provide flags
\code{'jacobian\_v'} or \code{'linop\_v'} respectively. Alternatively
a function handle to a function defined as \code{function res =
 evalFun(t, y)} can be provided instead of \code{'on'} to evaluate
the required parts in its own routine.

To add more flexibility we also enabled the user to provide his or her
own implementation to compute the matrix functions. This is especially
interesting for situations where the matrix functions can be computed
more cheaply, more easily, or in a structure-preserving way 
 due to special properties of the matrix. Then, 
instead of the \code{'direct'} or \code{'arnoldi'} settings a function
handle has to be given
\begin{sniplet}{}%
	function [ hOut, varargout ] = ... \\%
	~~~~~~~~matFun(job, t, y, h, flag, v, reusable, reuse, facs). \\%
\end{sniplet}
We refer to the documentation contained in the package for detailed
instruction. Here we will only give a short introduction.

Due to the possible complexity of this task, a number of stages have
to be incorporated into the process: in addition to the evaluation
itself there is an \code{'init'}, a \code{'registerjobs'} and an
\code{'initstep'} phase to precompute data -- indicated by the
corresponding value for \code{flag}.  It might happen that the
evaluator requires to reduce the step size to guarantee the prescribed
accuracy. Therefore it is possible to return a different value for
\code{hOut} than the input step size \code{h} to indicate such a time
step reduction.

Details of the implementation can be found in the manual. We suggest taking
the two \EXPODE internal evaluators -- \path{matFun/matFunDirect.m}
and \path{matFun/matFunKrylov.m} -- as a guideline.


\subsection{Custom Integrators}

The \EXPODE package was implemented such that it is open for user
specified extensions in many different ways. An example was pointed
out in the previous section for the evaluation of the matrix
functions. It is also possible to add new time integrators to
the package. To do so, one should use the developer package of
\EXPODE, which contains some useful tools for this purpose. The
integration steps for the integrators are written as plugins to the
\code{expode} routine. The developer does not need to worry about
things like direct user interaction syntactical options checking,
output control and other management tasks.

We provide scripts to generate an integrator \code{stub} and for
deployment. For a rudimentary integrator only three files have to be
edited after running the generation script: a setup routine which
gives some information to the \code{expode}, a routine that provides
the options for the user and the integration step itself. Detailed
information on the development process are available in the
documentation contained in the \EXPODE package.


\section{Examples}

In this section we want to discuss some numerical examples.

\subsection{van der Pol equation}

\begin{figure}[tbh]
	\includegraphics[width=0.48\textwidth]{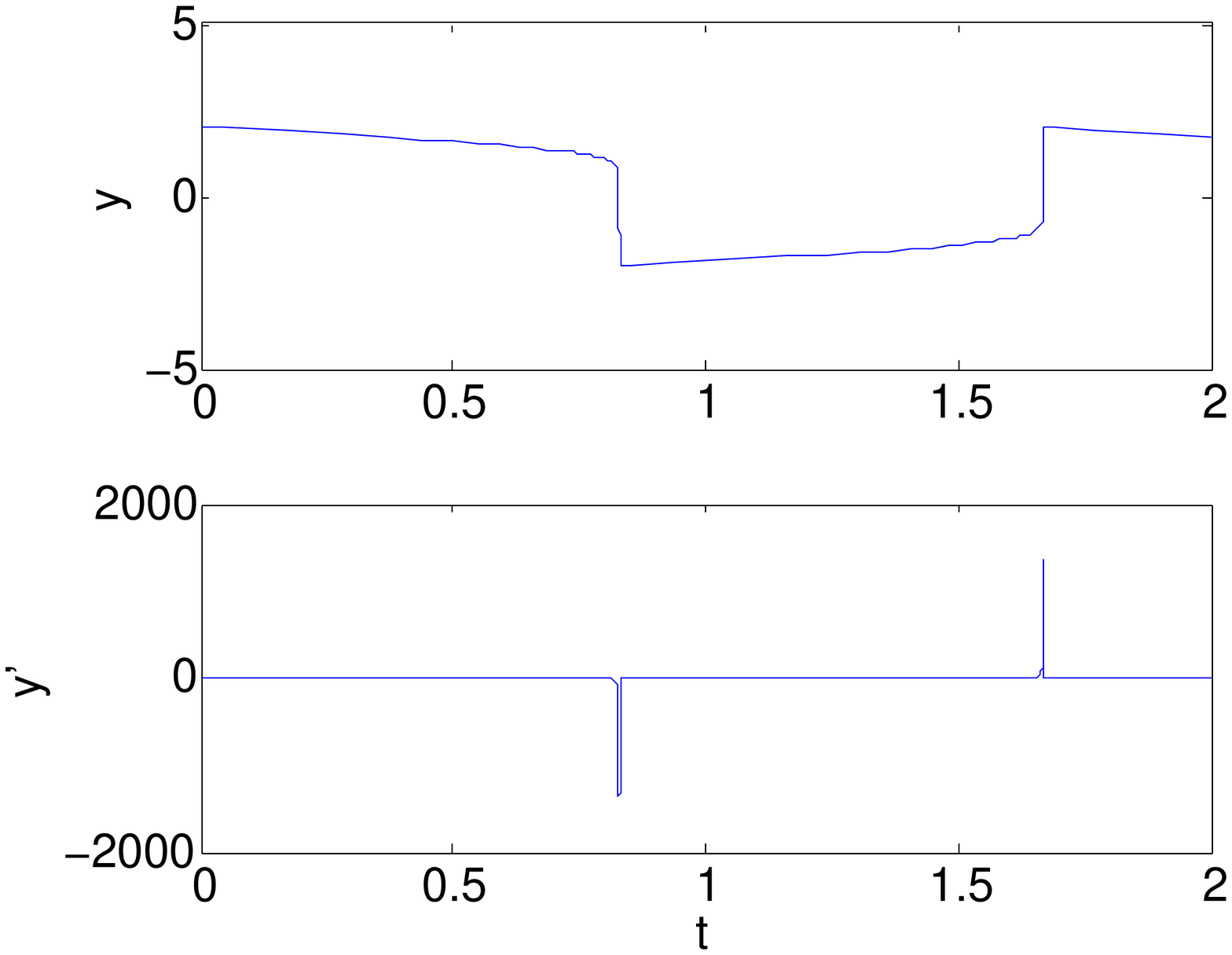}
	\quad
	\includegraphics[width=0.48\textwidth]{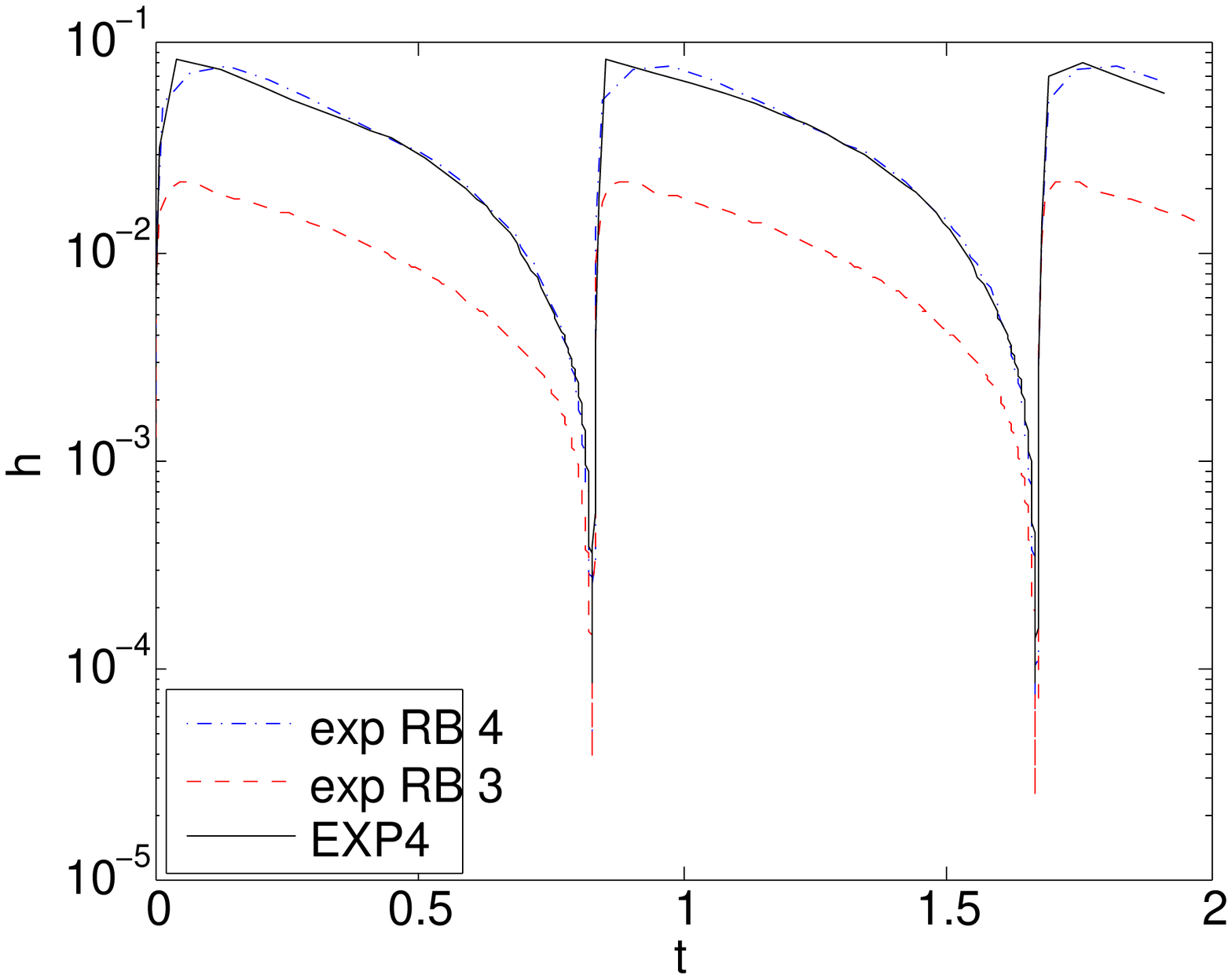}
	\caption{Left: $t$ vs. the solution $(y(t), y'(t))$ of the van der Pol equation ($\mu = 1000$), right: $t$ vs. step sizes $h$ chosen by different \EXPODE integrators with adaptive step size.}
	\label{fig:vdpSolution}
\end{figure}

To benchmark the adaptive step size implementation, let us consider
the well-known van der Pol equation \cite{vdP1920}. It describes a
non-linearly damped oscillator. Written as a second order \ode{} it reads
\[
	\uu''(t) - \mu(1 - \uu(t)^2) \uu'(t) + \uu(t) = 0.
\]
With large $\mu$ the equation becomes very stiff. A plot of the
solution with $\uu(0) = [2, -0.6]^T$ and $\mu = 1000$ can be found in
figure \ref{fig:vdpSolution}. At the vertical edges the first
derivative of the solution becomes very large so that small step sizes
have to be chosen. We plotted the step sizes selected
by the step size control of \EXP{rb4}, \EXP{rb3} and \EXP{4}. Due to
their larger (classical) order, \EXP{rb4} and \EXP{4} can choose
larger step sizes.

\subsection{Semilinear Problem}

\begin{figure}[tbh]
	\includegraphics[width=0.3\textwidth]{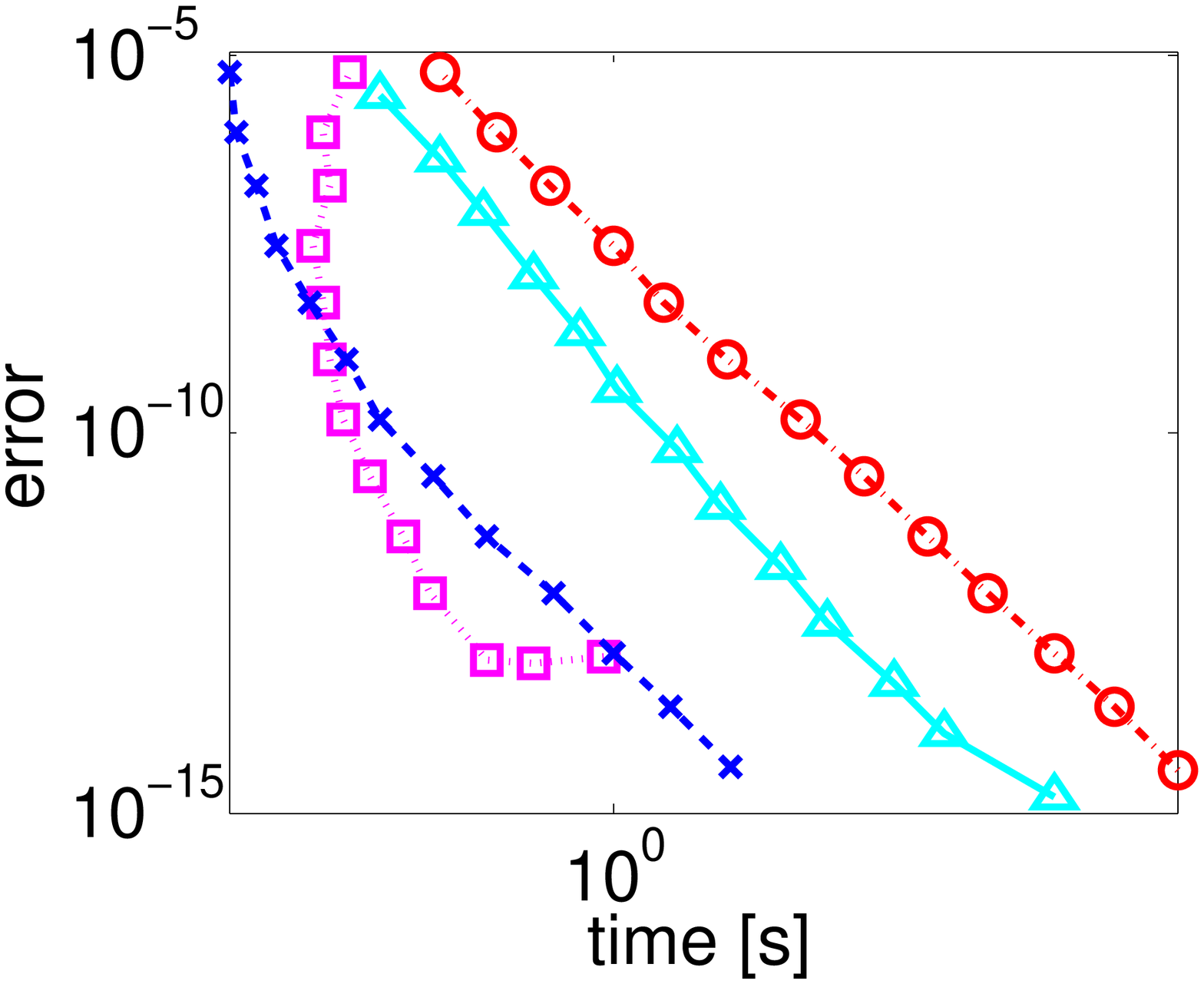}
	\quad
	\includegraphics[width=0.3\textwidth]{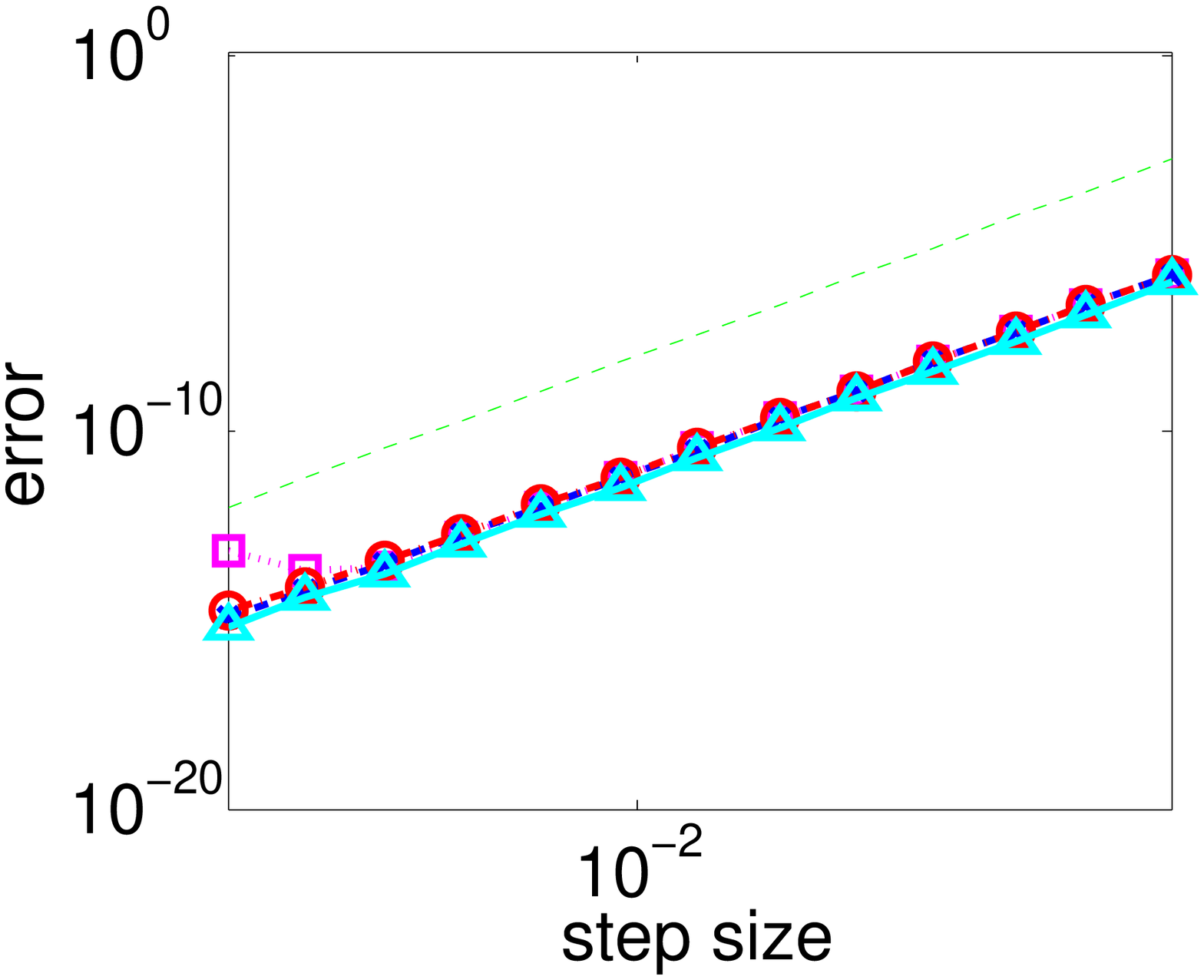}
	\quad
	\includegraphics[width=0.3\textwidth]{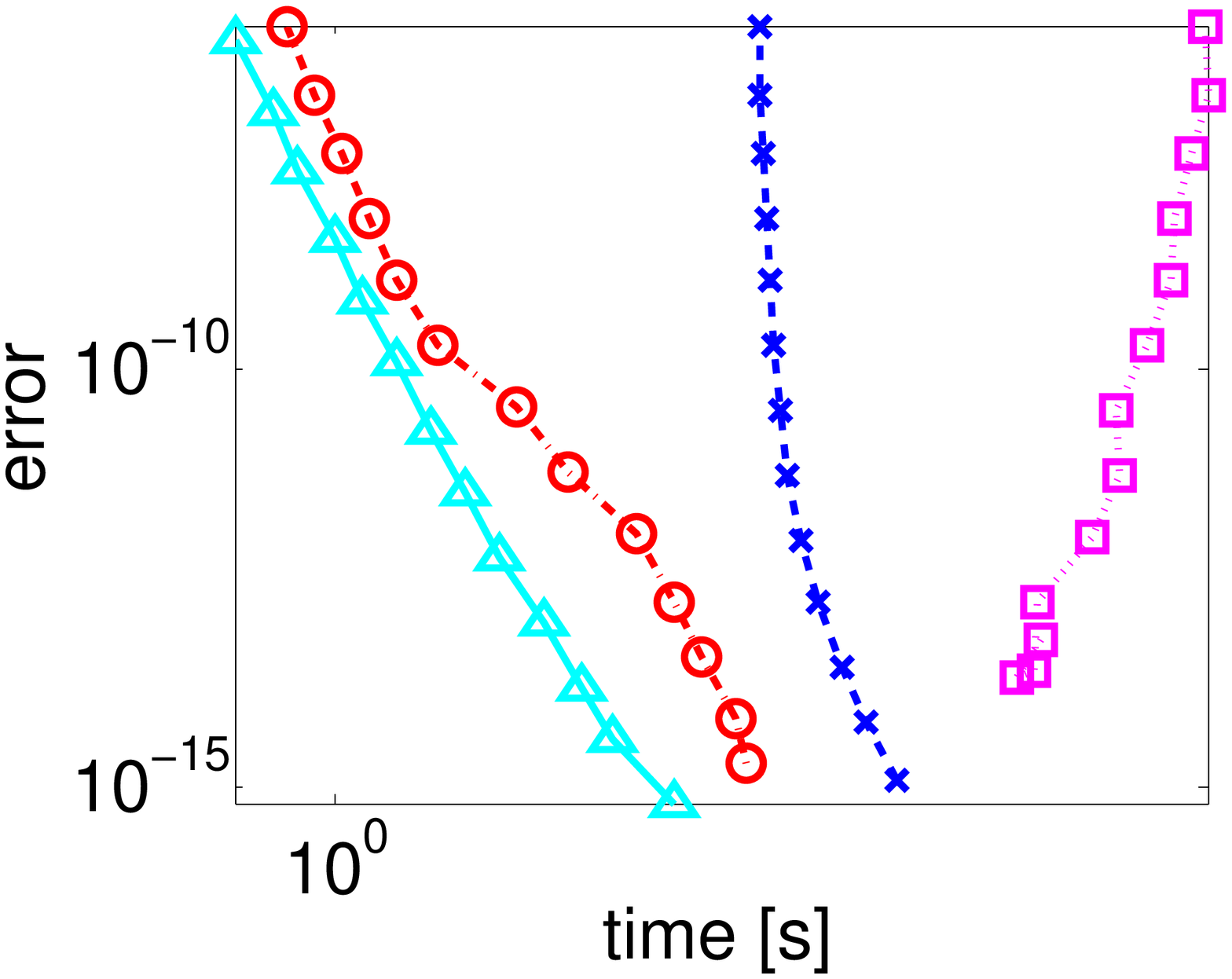}
	\caption{Left: Time vs.\ error plots for the one-dimensional version, middle: step size vs.\ error and right: time vs.\ error plots for the two-dimensional version of the semilinear problem \eqref{eq:semilinExample}. Dotted: \EXPINT, dashed: \EXPODE with diagonalization, dash-dotted: \EXPODE with Krylov method for the Matrix functions, solid: \EXPRB, thin dashed in the middle image: slope line of order four.}
	\label{fig:semi1Results}
\end{figure}

Another \MATLAB package implementing exponential 
Runge-Kutta methods and exponential  multistep methods, is \EXPINT
\cite{Ber2007}. Note that the evaluation of matrix functions in
\EXPINT corresponds to the ``direct'' option in our package. This
favors our package for large scale problems.
We compare the efficiency of both packages on the semilinear problem
\begin{equation}
	\label{eq:semilinExample}
	\uu' = \Delta \uu + \frac 1 {1 + \uu^2} + \Phi,
	\quad
	\uu(0, x, y) = x (1 - x) y (1 - y)
\end{equation}
with homogeneous Dirichlet conditions on $\Omega = [0, 1]^2$. $\Phi$
is chosen such that the exact solution is given by $\uu(t, x, y) = x
(1 - x) y (1 - y) \exp(t)$, see e.g. \cite{Ber2007} or
\cite{HocO05SIAM}. The space discretization was done with
finite differences with $N = 50$ inner grid points in each
dimension. This yields a system of \ode{}s of medium size dimension $N^2 = 2500$.
We compared \EXPINT's and our implementation of the
Hochbruck-Ostermann exponential Runge-Kutta scheme \cite{HocO05SIAM}
and our implementation of the three stage exponential Rosenbrock
method. 
Additionally we solved the one-dimensional version of this problem with $N = 100$ grid points to retrieve a low-dimensional \ode{} of $N = 100$ degrees of freedom.
The results are shown in figure \ref{fig:semi1Results}.

In the left image we see that for low dimensional systems and small target accurcies, \EXPINT is clearly preferable, since it is much faster. Especially the Krylov matrix function evaluator has negative impact on the performance since it has the overhead of building up the Krylov subspaces each time step while \EXPINT only computes some matrix exponentials via Pad\'e approximation. For $100 \times 100$ matrices the solution of the linear systems arising there is quite cheap.

In the middle image we see numerical order four of all solvers.

\EXPINT uses Pad\'e approximations and a scaling and squaring technique to compute the $\varphi$ functions \eqref{eq:phi-functions}. Let $\AA$ be the discretized Laplacian operator, then $h \AA$ has to be scaled, such that its norm is smaller than a given bound. This implies that smaller time step sizes lead to fewer scaling and squaring steps. After computing the $\varphi$ functions of the matrix once only matrix vector products are needed for the time stepping, such that a smaller time step size actually leads to faster computations.

For the direct solver in \EXPODE we have an initialization phase, where we diagonalize the operator $h A$, which takes the same time independent of $h$. Then we apply $\varphi_k(h A)$ directly as \EXPINT does, so we get only slightly increased costs for higher accuracy by the smaller time step size.

In \EXPODE with the Arnoldi method for the matrix functions we build
up Krylov subspaces in each timestep and cannot reuse data from
previous time steps. Nevertheless, this is much faster than computing
the matrix functions themselves.

\subsection{Brusselator}

As our next example we use the Brusselator example, which was also
used to benchmark the original \EXP4 code \cite{HocLS98}. It models a
reaction-diffusion process with two species $u = u(t, x, y)$ and $v =
v(t, x, y)$. The
equation reads
\begin{align*}
    \frac{\dd}{\dd t} \left( \begin{matrix} u \\ v \end{matrix} \right) &
    = A \left( \begin{matrix} u \\ v\end{matrix} \right)
        + \left( \begin{matrix} \gamma + u^2 v u \\ u - u^2
            v \end{matrix}   \right),
    \qquad A = \left[ \begin{matrix} \alpha \Delta - (\beta + 1) & 0 \\ \beta & \alpha \Delta \end{matrix} \right]
\end{align*}
together with homogeneous Neumann boundary conditions. The initial
value for $u$ is taken from \MATLAB's \code{peaks} function and $v(0, \cdot,
\cdot) \equiv 0$. The stiffness results from the diffusion term $A$. We choose the parameters $\gamma = 1$, $\beta = 3.4$ and the diffusion coefficient $\alpha = 10^{-2}$. Space discretization is done via finite differences with $100^2$ grid points in $\Omega = [0, 1]^2$. This results in a system of 20.000 unknowns.

\begin{figure}[tbh]
	\includegraphics[width=0.48\textwidth]{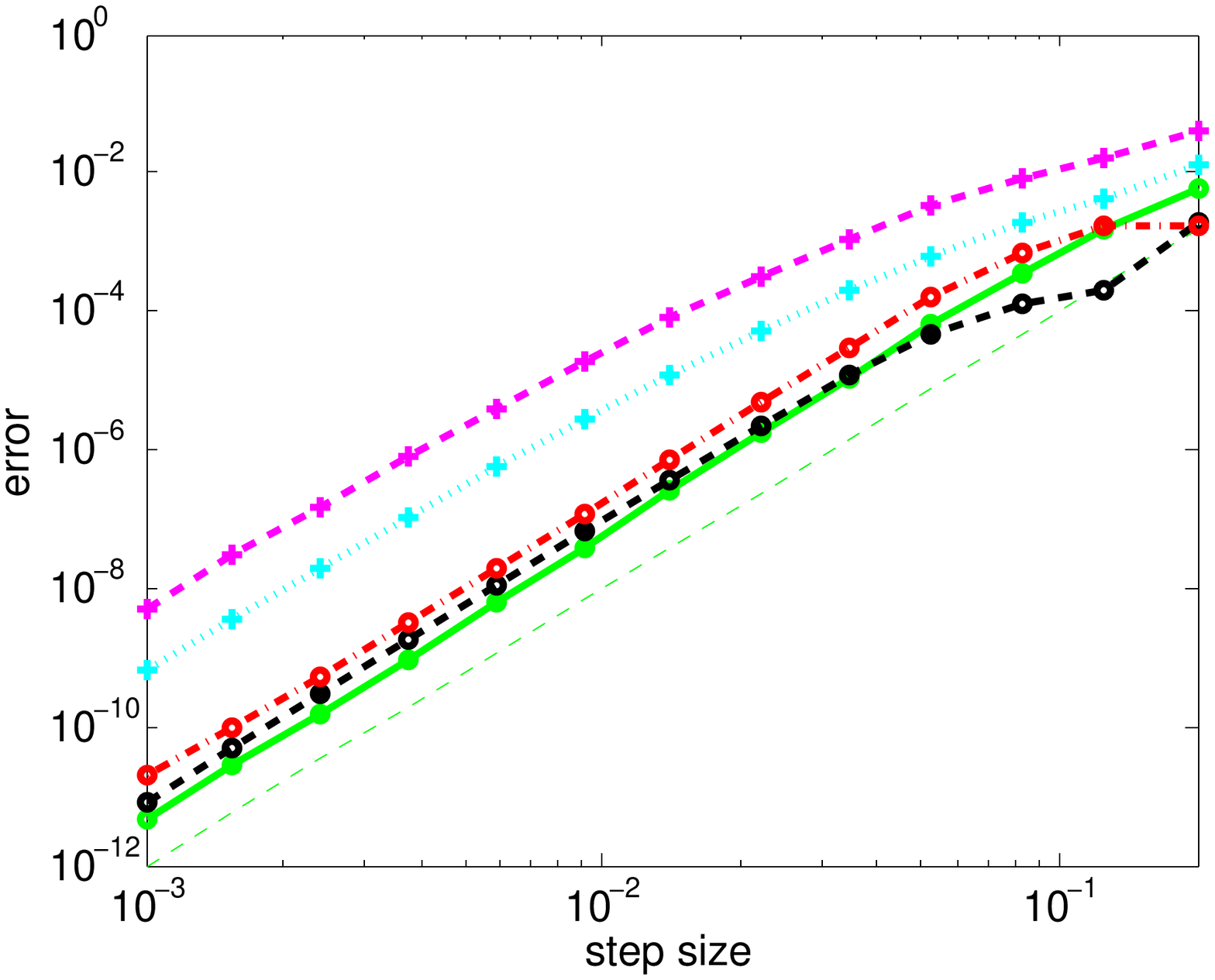}
	\quad
	\includegraphics[width=0.48\textwidth]{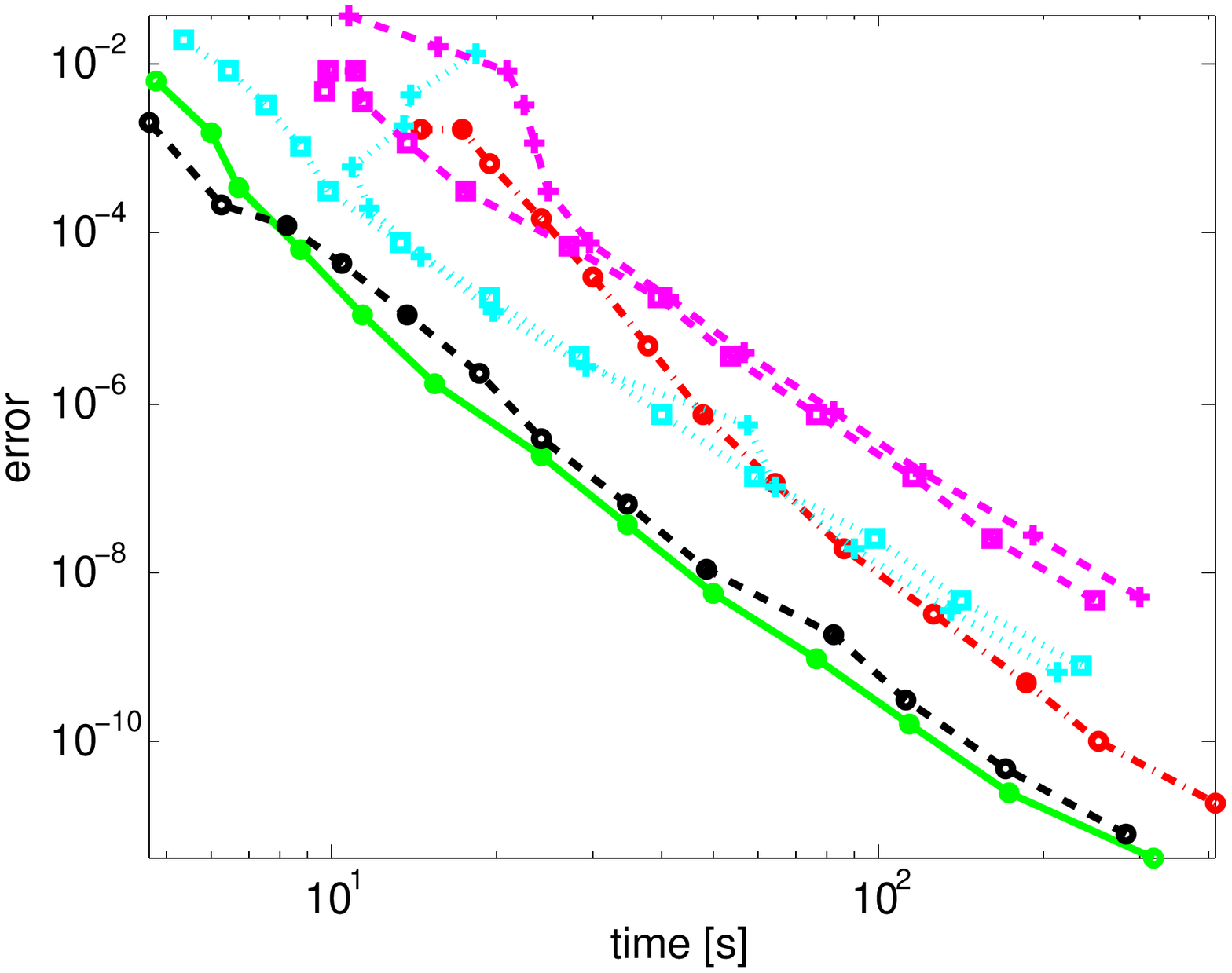}
	\caption{Step size vs.\ error and time vs.\ error plots for
          different fourth order integrators applied to the
          Brusselator example. Integrators: 3 stage \EXPRB (solid
          line), \EXP4 (dashed line), 5 stage Hochbruck-Ostermann
          exponential Runge-Kutta-Scheme (dash-dotted line), 3 step
          \EXPMS with  initial values computed by fixed point
          iteration (dotted, plus markers) and with \EXPRB initial values (dotted, square markers) and 4 step \EXPMSSEMI with Runge-Kutta initial values (dashed, plus markers) and with initial values computed by fixed point
          iteration (dashed, square markers).}
	\label{fig:BrussExample}
\end{figure}

In figure~\ref{fig:BrussExample} we present time step vs.\ error and
work-precision diagrams for the fourth order representatives of each of
our integrator classes. The multistep integrators use the fixed point
iteration for their starting values. As matrix function evaluator we
choose the Arnoldi Krylov method with a maximal dimension of  36 for
the Krylov subspaces
for the linearized one-step methods and 100 for the others. We choose
to use constant step sizes with matching tolerances for the Arnoldi
process, since not all integrators allow adaptive step sizes. The
error is measured in the maximum norm against a reference solution
computed by \MATLAB's {\tt ode15s} with sufficiently high accuracy
requirements. This means that we only consider the \ode{} error and not 
the spatial error. We added dashed lines for slope of order four.

In the time vs.\ error plot we arranged data for the same runs and
additionally gave two curves for the multistep integrators using
\EXPRB (linearized) and \EXPRK (exponential Adams) for the initial
values.

\begin{figure}[tbh]
	\includegraphics[width=0.48\textwidth]{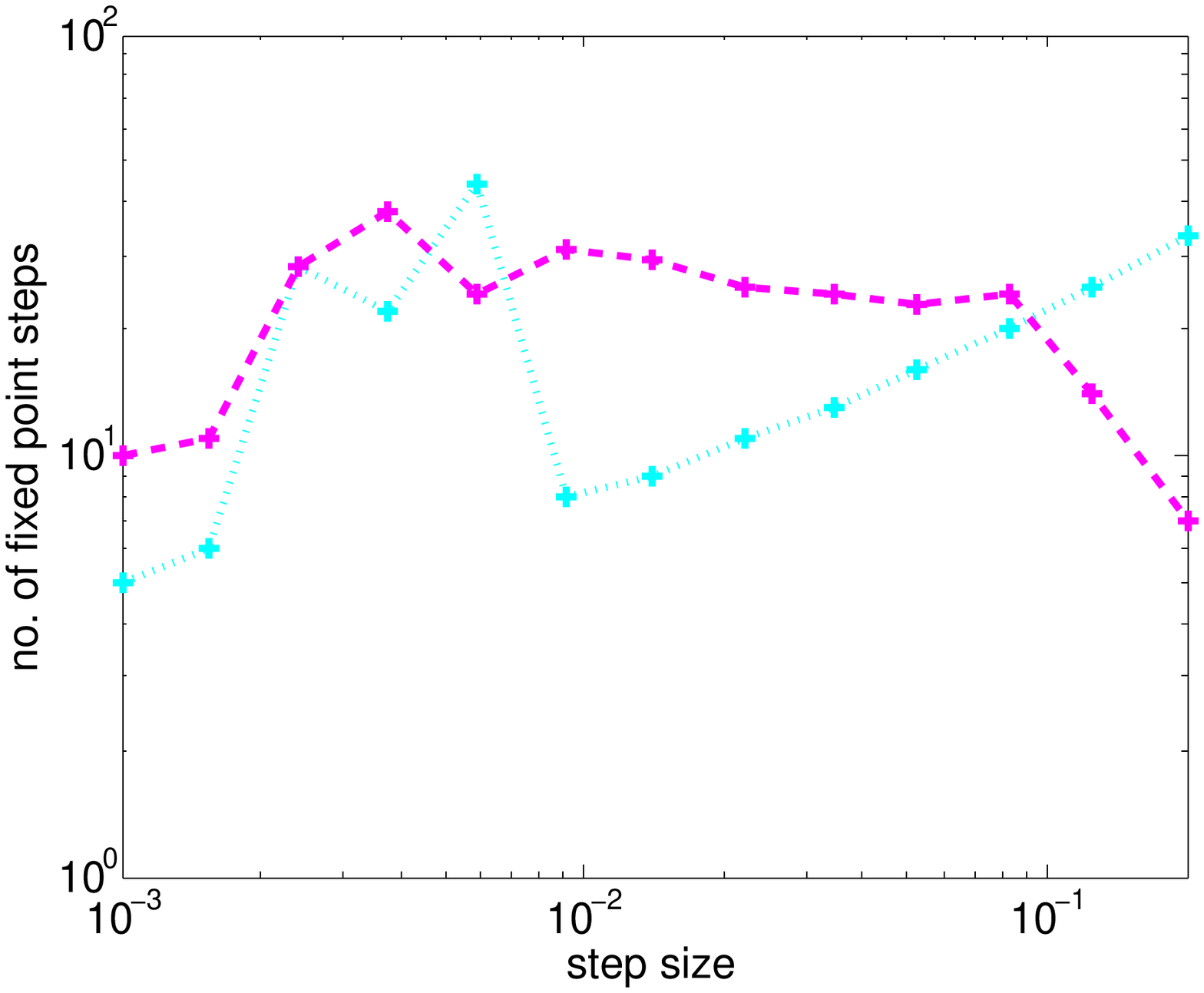}
	\quad
	\includegraphics[width=0.48\textwidth]{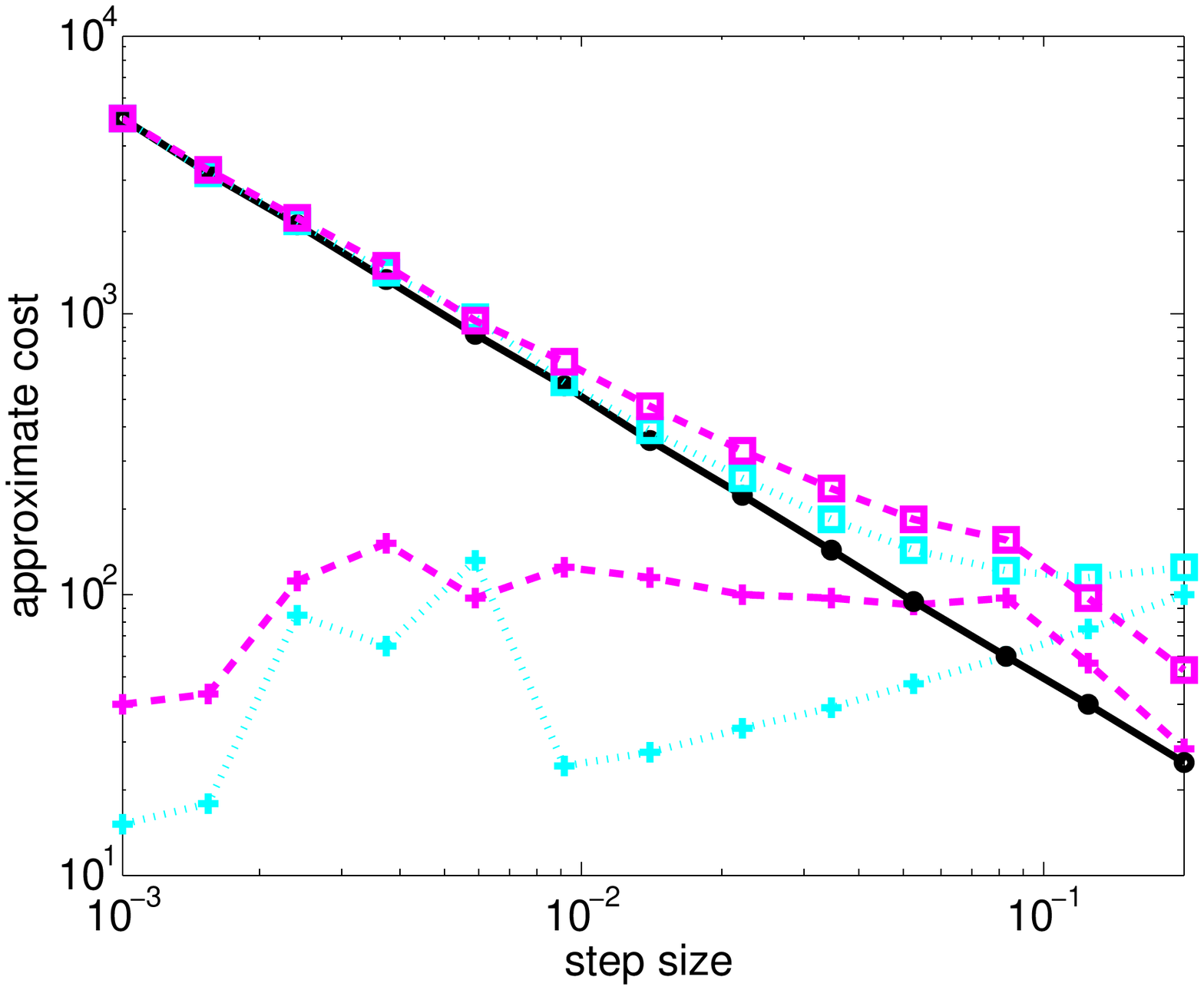}
	\caption{Left: number of fixed point steps required to fulfill the accuracy requirement for the different step sizes in the Brusselator Example. \EXPMS: dotted, plus markers, \EXPMSSEMI: dashed, plus markers. Right: approximate costs for the fixed point iterator linestyles and markers as in left picture, approximate costs for the (non-startup) time steps, approximate costs for the whole integration, line styles as left, now square markers.}
	\label{fig:BrussFixedpoint}
\end{figure}

In the first image we observe fourth order convergence for all schemes as $h \to 0$. The two linearized one-step methods have the best error constants (offset on the error-axis) closely followed by the exponential Runge-Kutta scheme.

In the second one we see the linearized one step methods lead to the best
accuracy for a given computation time. Here the exponential
Runge-Kutta scheme lags a bit further behind, since it is
computationally more expensive due to its five internal stages. \EXP4
has seven of them, but the first three and the second three of them can be computed simultaneously, see \cite{HocLS98}, so \EXP{rb4} and \EXP{4} have
approximately the same computational costs. The multistep integrators
(with plus sign markers) are much weaker for larger time
steps. 
The reason for this is the fixed point iteration which gets relatively
more expensive when the time step size increases. In figure
\ref{fig:BrussFixedpoint} we show the number of fixed point steps
required versus the step size chosen and an approximation of the
normal time step costs versus the fixed point step costs. Each fixed point
step costs approximately as much as the number of backwards steps in use
times the cost of one time step. So the \EXPMS solver fixed point step
costs 3 normal time steps, the \EXPMSSEMI fixed point step 4 normal time
steps. The plot assumes constant time for the matrix function
evaluator.  We also gave curves for multistep methods with
\EXPRB/\EXPRK initial values, which are significantly more efficient,
especially for larger step sizes in this case.
Note that multistep methods of order higher than four can only be obtained using the fixed point iteration.

\subsection{2D Maxwell's equations with spatial discontinuos Galerkin
  discretization}

\def\ratioA{$1$}
\def\avgTimeStepQuotA{8.06053}
\def\avgSpeedupA{1.23683}

\def\ratioB{$\left( \frac{1}{4} \right)^2$}
\def\avgTimeStepQuotB{15.481}
\def\avgSpeedupB{2.23478}

\def\ratioC{$\left( \frac{1}{16} \right)^2$}
\def\avgTimeStepQuotC{26.1532}
\def\avgSpeedupC{3.39336}

\def\ratioD{$\left( \frac{1}{64} \right)^2$}
\def\avgTimeStepQuotD{35.5925}
\def\avgSpeedupD{4.46899}

\def\ratioE{$\left( \frac{1}{256} \right)^2$}
\def\avgTimeStepQuotE{44.667}
\def\avgSpeedupE{4.94512}

As an example for a hyperbolic problem we consider Maxwell's equations
in two space dimensions. For a charge free domain the
equations are given by
\begin{align*}
	\mu \frac{\partial}{\partial t} \vec H
	& = - \nabla \times \vec E \\
	\epsilon \frac{\partial}{\partial t} \vec E
	& = \nabla \times \vec H.
\end{align*}
Here $\vec E$ denotes the electric field, $\vec H$ the magnetic field,
$\epsilon$ the permittivity and $\mu$ the permeability. Assuming
constant $\mu$ and $\epsilon$
and using normalized cartesian coordinates we can
eliminate the two material parameters.

\begin{figure}[tbh]
	\centering
    \includegraphics[width=0.35\textwidth]{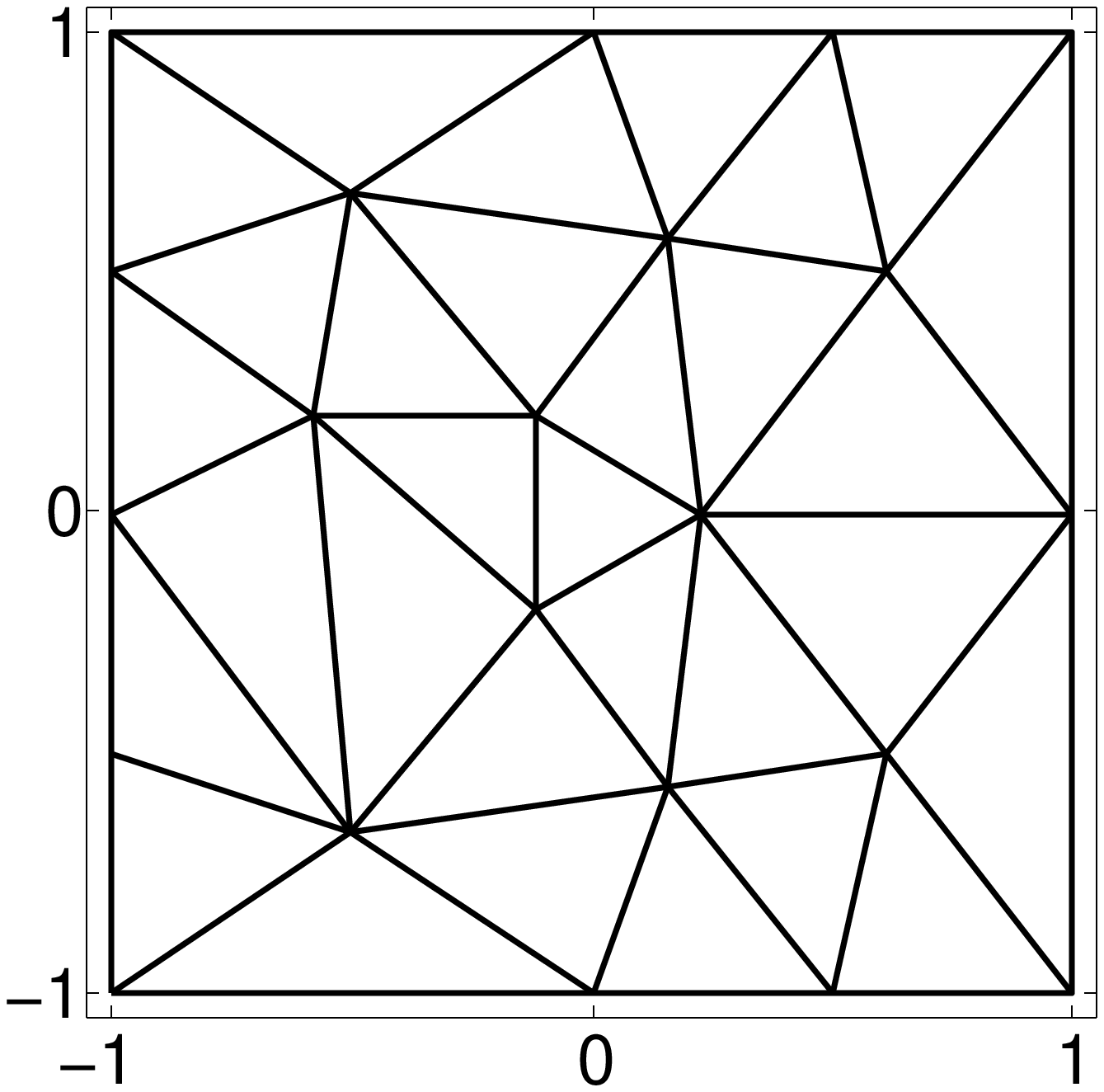}
    \quad
    \includegraphics[width=0.35\textwidth]{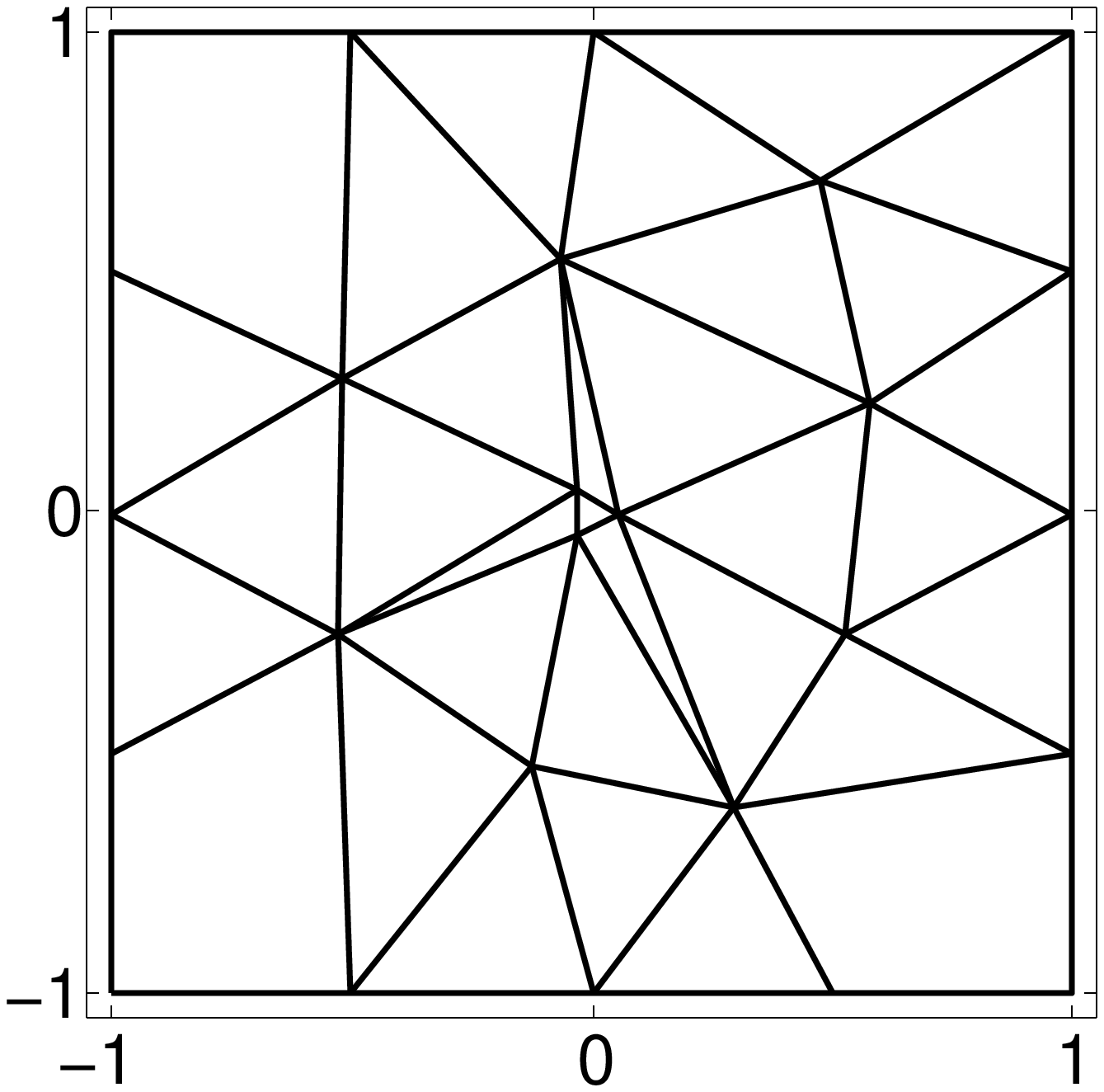}
    \caption{Examples for the grids used in the Maxwell example with ratios \ratioA{} (left) and \ratioB{} (right).}
    \label{fig:MaxwellGridExample}
\end{figure}

As domain we use a magnetic box $\Omega = [-1, 1]^3$.
Our initial conditions are chosen, such that 
\begin{align*}
	\vec E & = \left[\begin{matrix}
			0 \\
			0 \\
			\sin(k \pi x) \cdot \sin(k \pi y) \cdot \cos(t k \sqrt 2 \pi)
	\end{matrix}\right] \\
	\vec H & = \left[\begin{matrix}
			\frac {-1}{\sqrt 2} \sin(k \pi x) \cdot \cos(k \pi y) \cdot \sin(t k \sqrt 2 \pi) \\
			\frac 1   {\sqrt 2} \cos(k \pi x) \cdot \sin(k \pi y) \cdot \sin(t k \sqrt 2 \pi) \\
			0
	\end{matrix} \right]
\end{align*}
is the exact solution. We assume the arising quantities to be constant
in the $z$-direction to reduce to two spatial dimensions. The space
discretization is done using the discontinuos Galerkin method (DG-method)
from the codes by \cite{hesthaven2008nodal}. We created grids with a tiny
triangle in the center to simulate very filigree structures and bad quality
grids. We run four different ratios between maximal area of the outer
elements to area of the center triangle, see figure \ref{fig:MaxwellGridExample}
for examples. For grid generation we use \cite{Shewchuk96}.

We were interested how the existence of tiny elements in an otherwise mostly
coarse grid would affect the stability requirements by an exponential integrator.

In our numerical experiment we created time vs.\ error plots. The error was measured in the $L^2$-Norm against the exact solution.
The spatial resolution was refined and the time step sizes were chosen such that the solvers are stable and the time and the spatial discretization errors are almost equal.

We compared the explicit space saving order four Runge-Kutta method used in \cite{hesthaven2008nodal} with our three stage exponential Rosenbrock-type (\EXPRB) and the \EXP4 integrators, both using the Arnoldi method to compute the matrix functions. The time step size in the Runge-Kutta solver is automatically chosen, such that the solver is stable.
This automatically leads to a time discretization error in the magnitude of the spatial error.
For \EXPRB and \EXP4 we set the accuracy requirement to the values retrieved by a reference solution and allowed adaptive step size choice. We also used the $L^2$ norm for the step size estimator. This was done providing a custom error norm function using the mass-scalarproduct from the DG-codes via the \refOption{NormFunction} option.

\begin{figure}[tbh]
	\centering
	\def\figsize{0.23\textwidth}

    \mbox{%
    	\vtop{\vskip0pt\hbox{%
    		\includegraphics[width=\figsize]{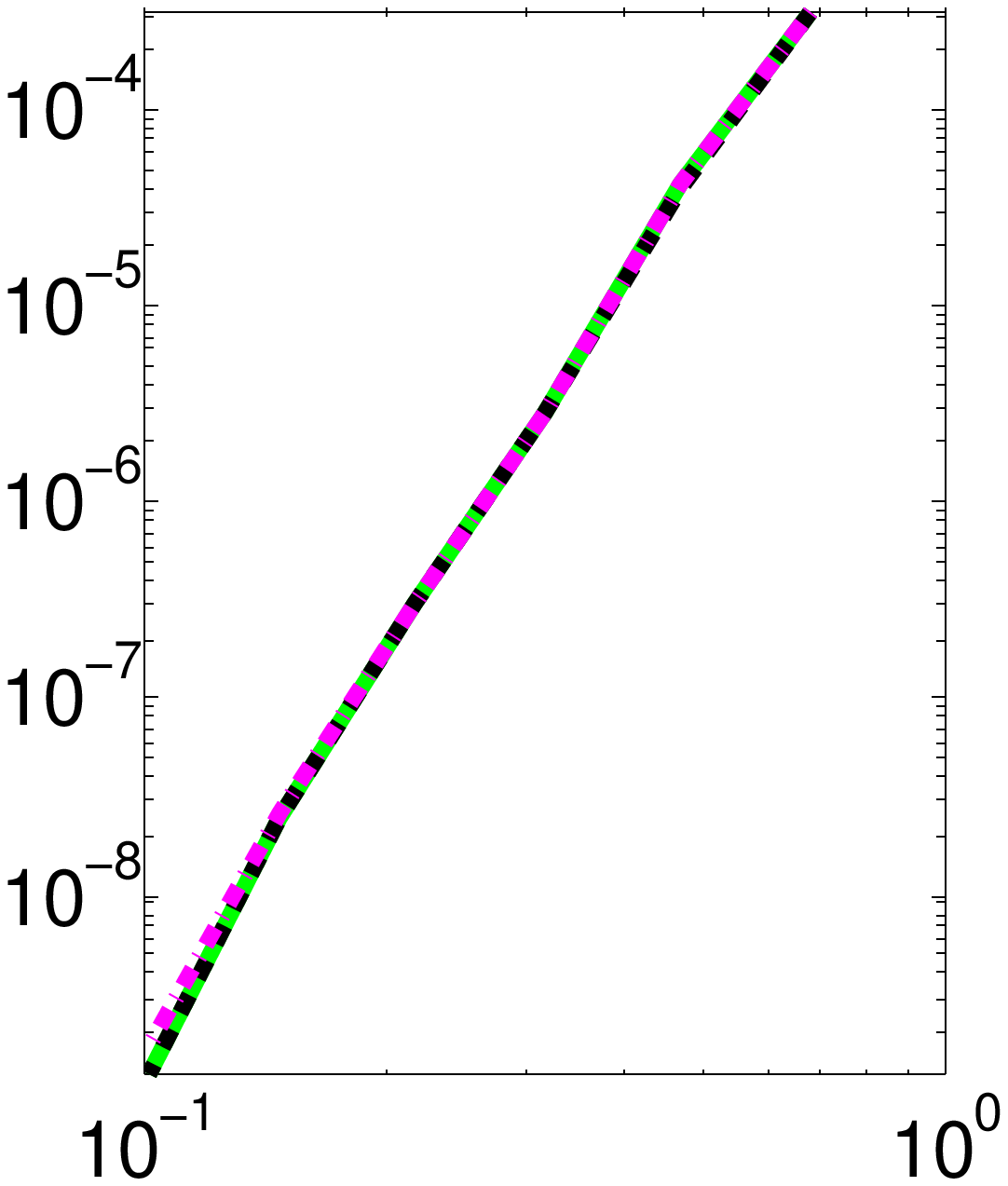}%
    	}}%
    	\;%
	    \vtop{\vskip0pt\hbox{%
	    	\includegraphics[width=\figsize]{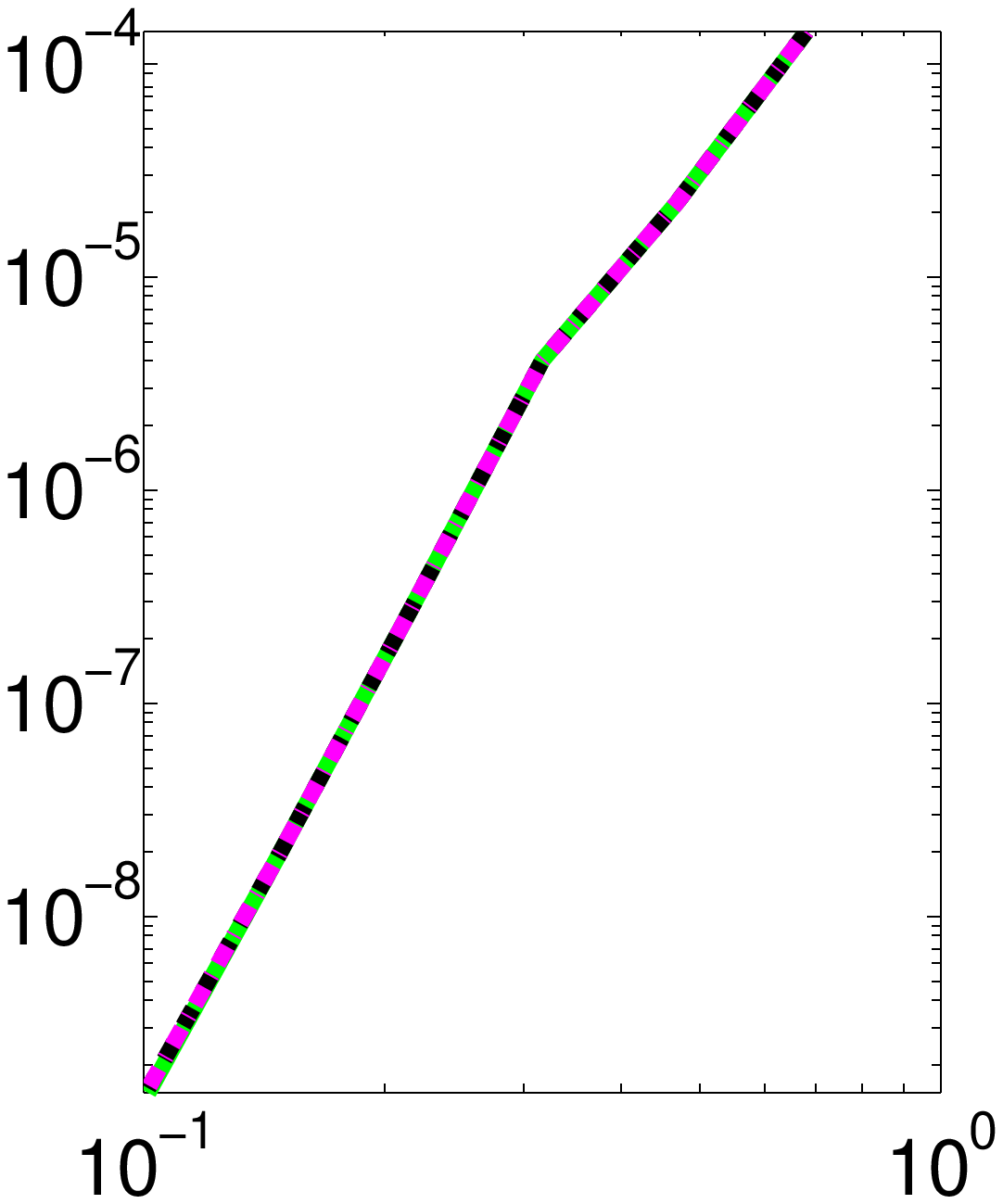}%
	    }}%
    	\quad%
		\vtop{\vskip0pt\hbox{%
			\includegraphics[width=\figsize]{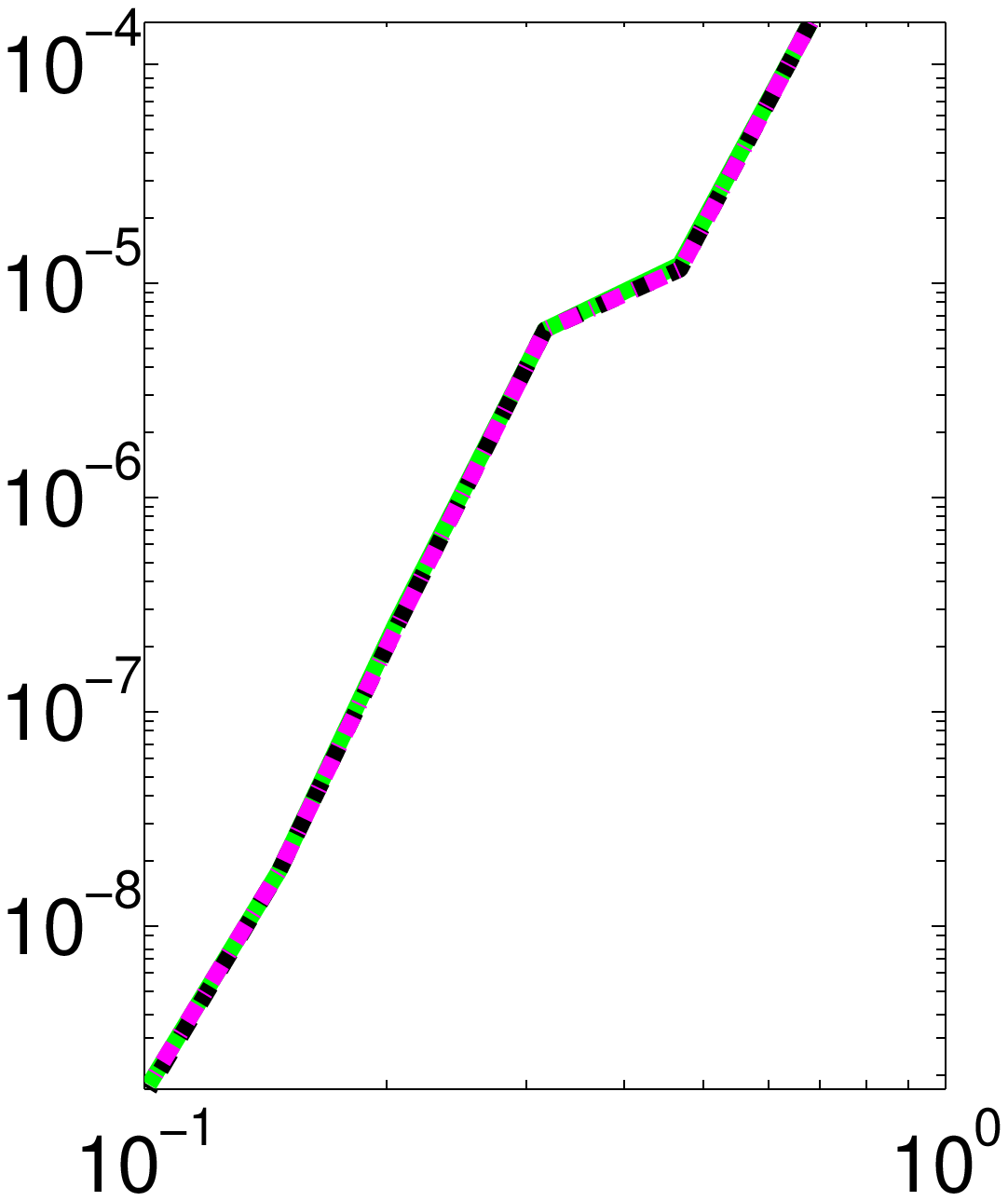}%
		}}%
    	\;%
    	\vtop{\vskip0pt\hbox{%
	    	\includegraphics[width=\figsize]{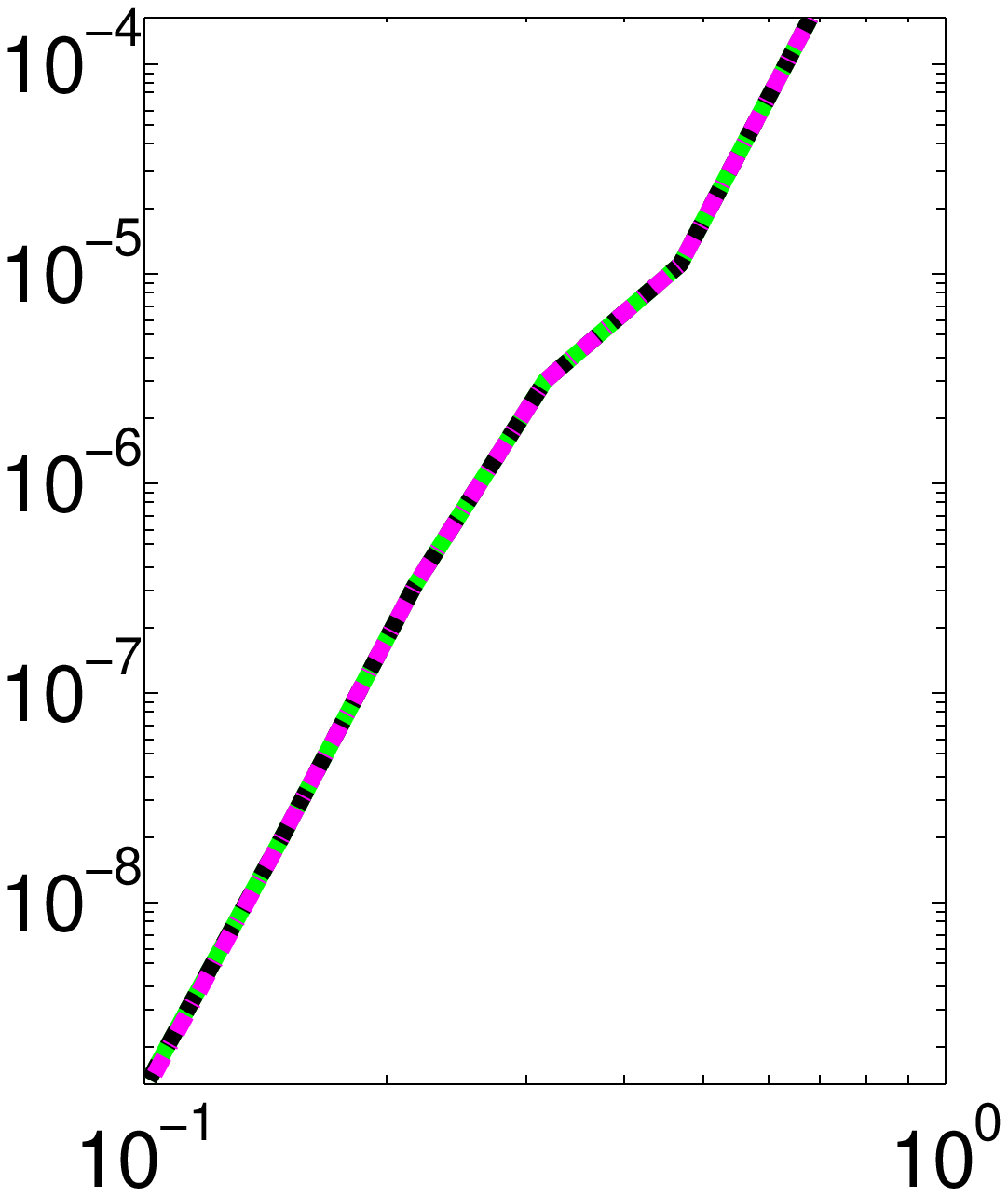}%
	    }}%
	}\\[0.25em]
	$\Delta x$ vs.\ error ($L^2$-Norm) \\[0.5em] 

    \mbox{%
    	\vtop{\vskip0pt\hbox{%
    		\includegraphics[width=\figsize]{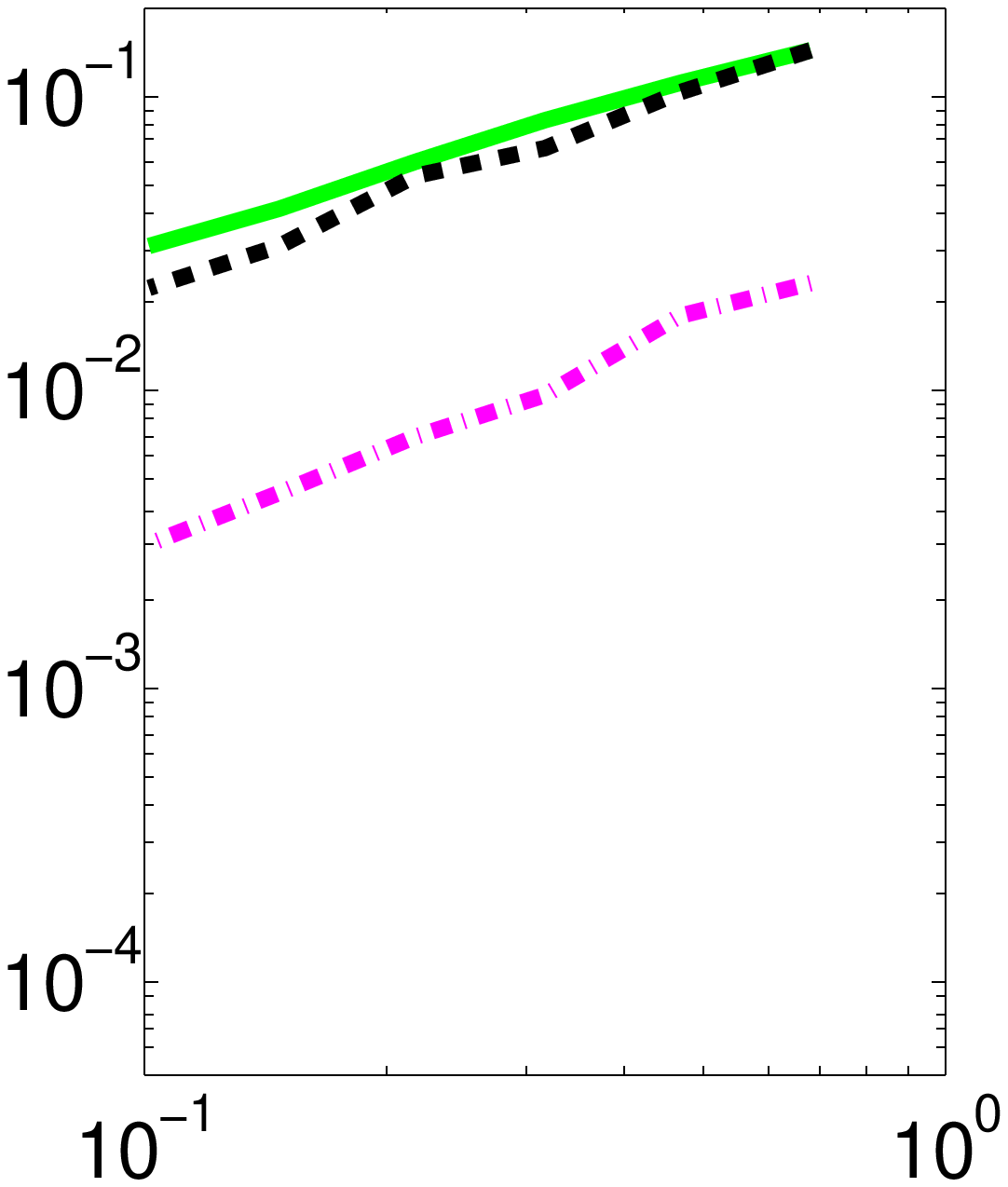}%
    	}}%
    	\;%
		\vtop{\vskip0pt\hbox{%
			\includegraphics[width=\figsize]{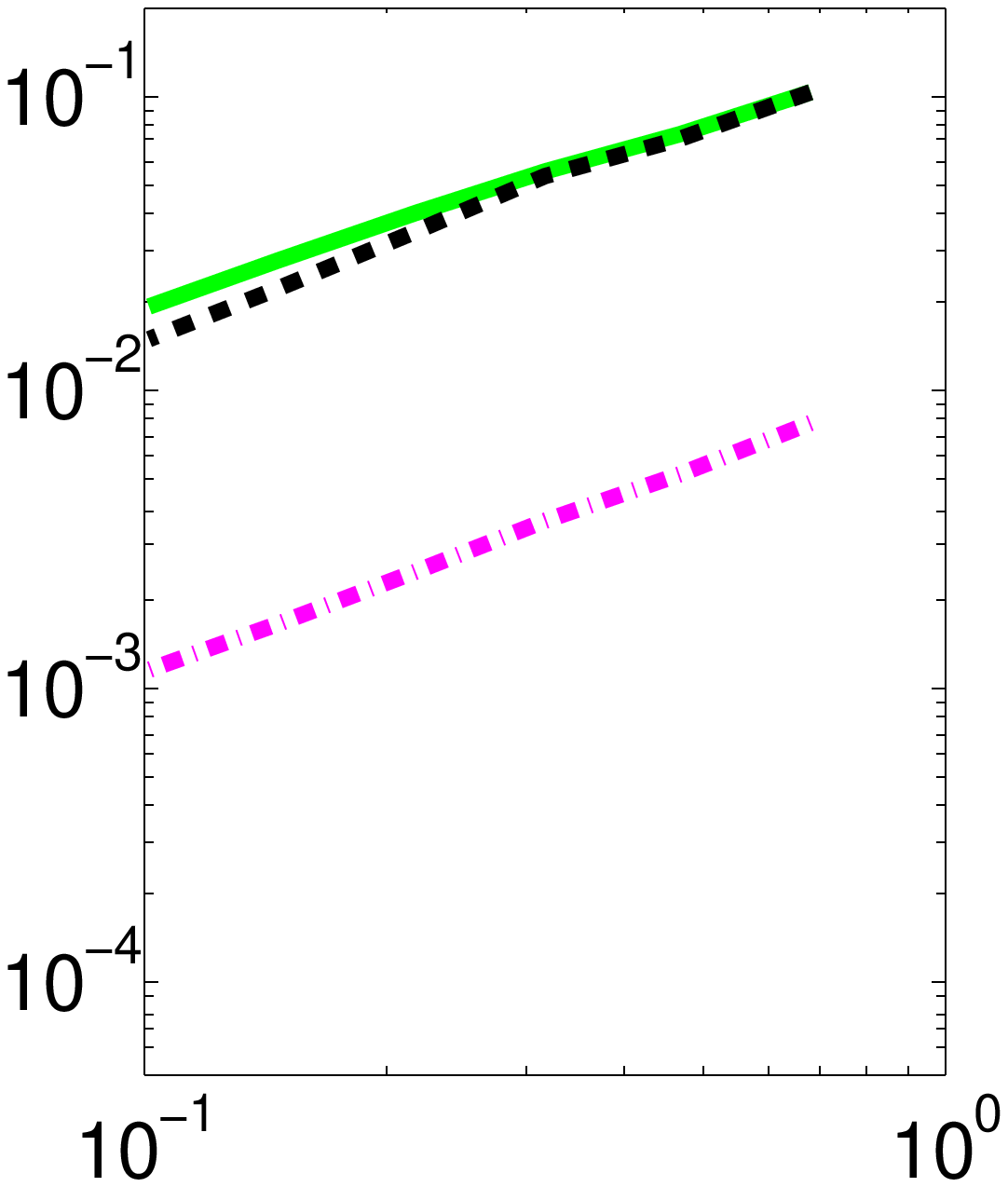}%
		}}%
    	\;%
    	\vtop{\vskip0pt\hbox{%
		    \includegraphics[width=\figsize]{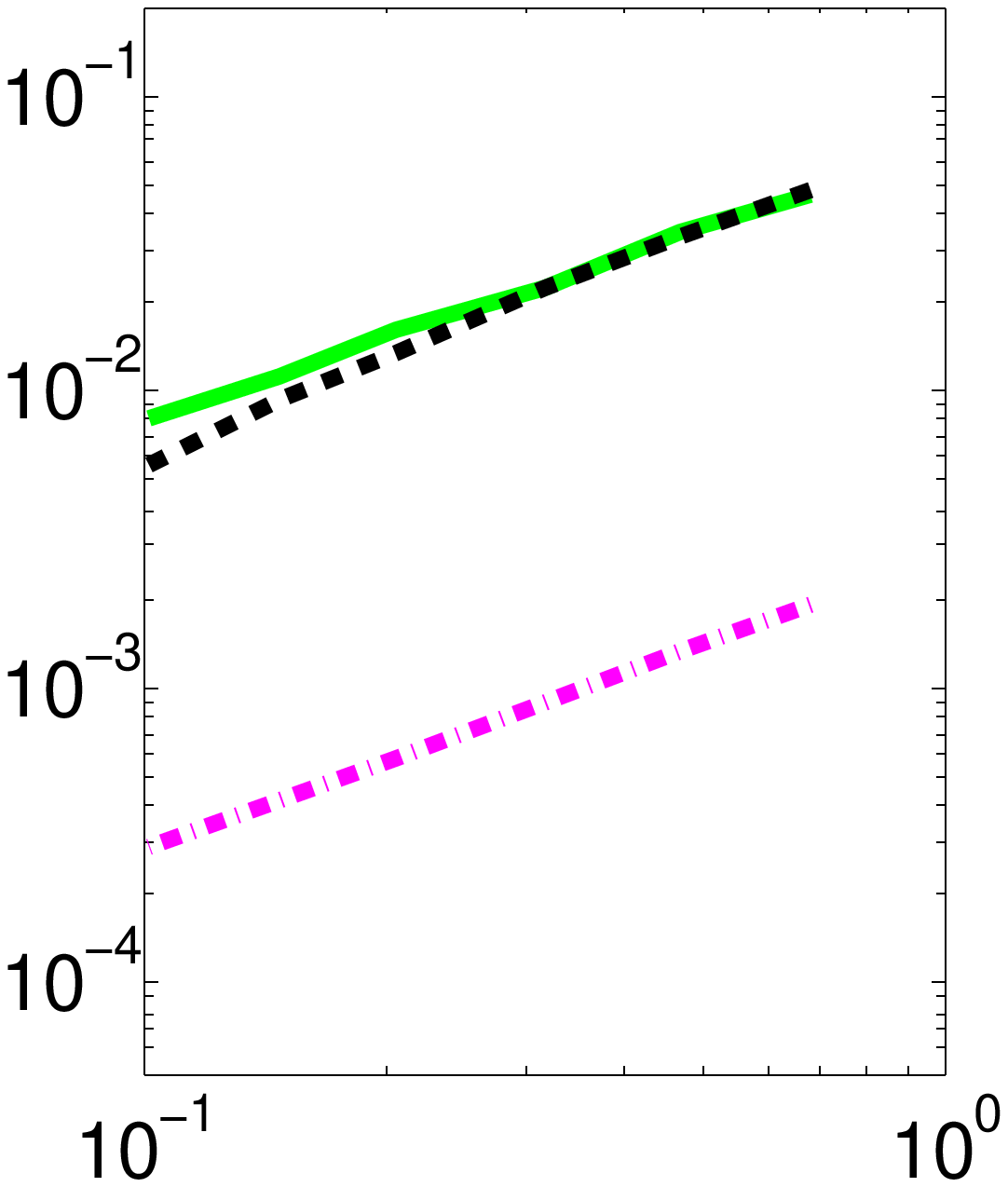}%
		}}%
    	\;%
		\vtop{\vskip0pt\hbox{%
	    	\includegraphics[width=\figsize]{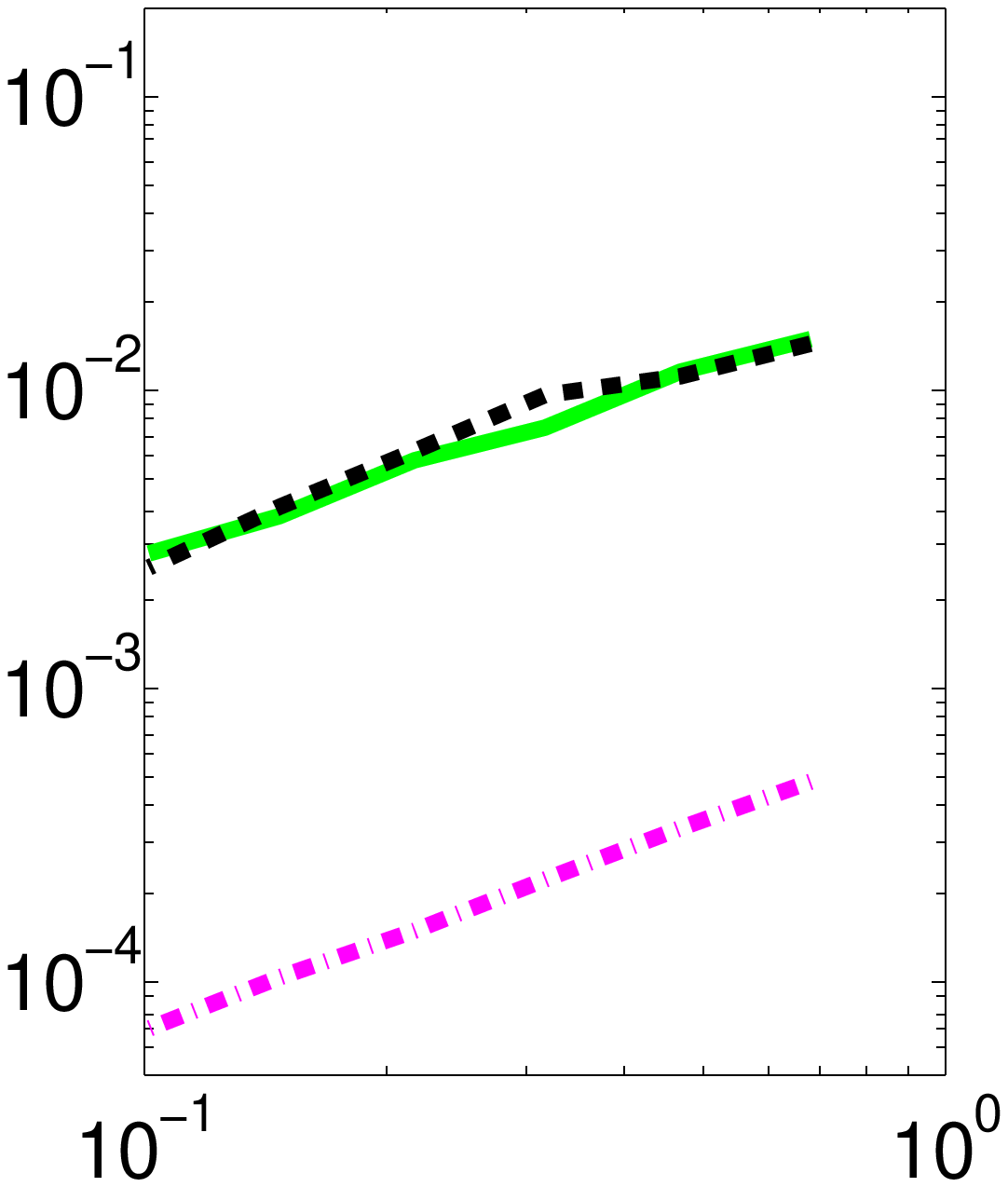}%
		}}%
	}\\[0.25em]
    $\Delta x$ vs.\ $h$ \\[0.5em]
    
    \mbox{%
	    \vtop{\vskip0pt\hbox{%
	    	\includegraphics[width=\figsize]{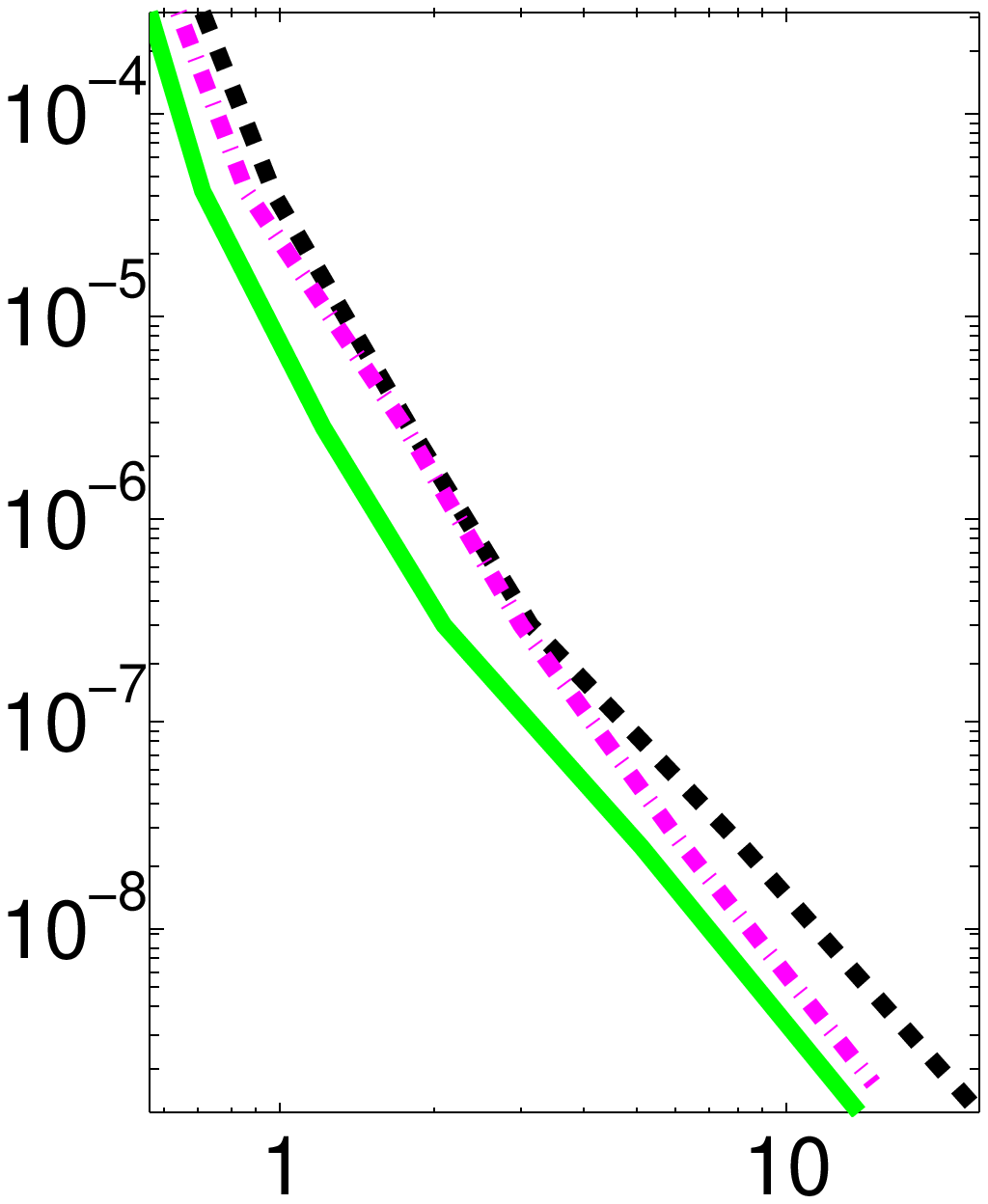}%
	    }}%
		\;%
    	\vtop{\vskip0pt\hbox{%
	    	\includegraphics[width=\figsize]{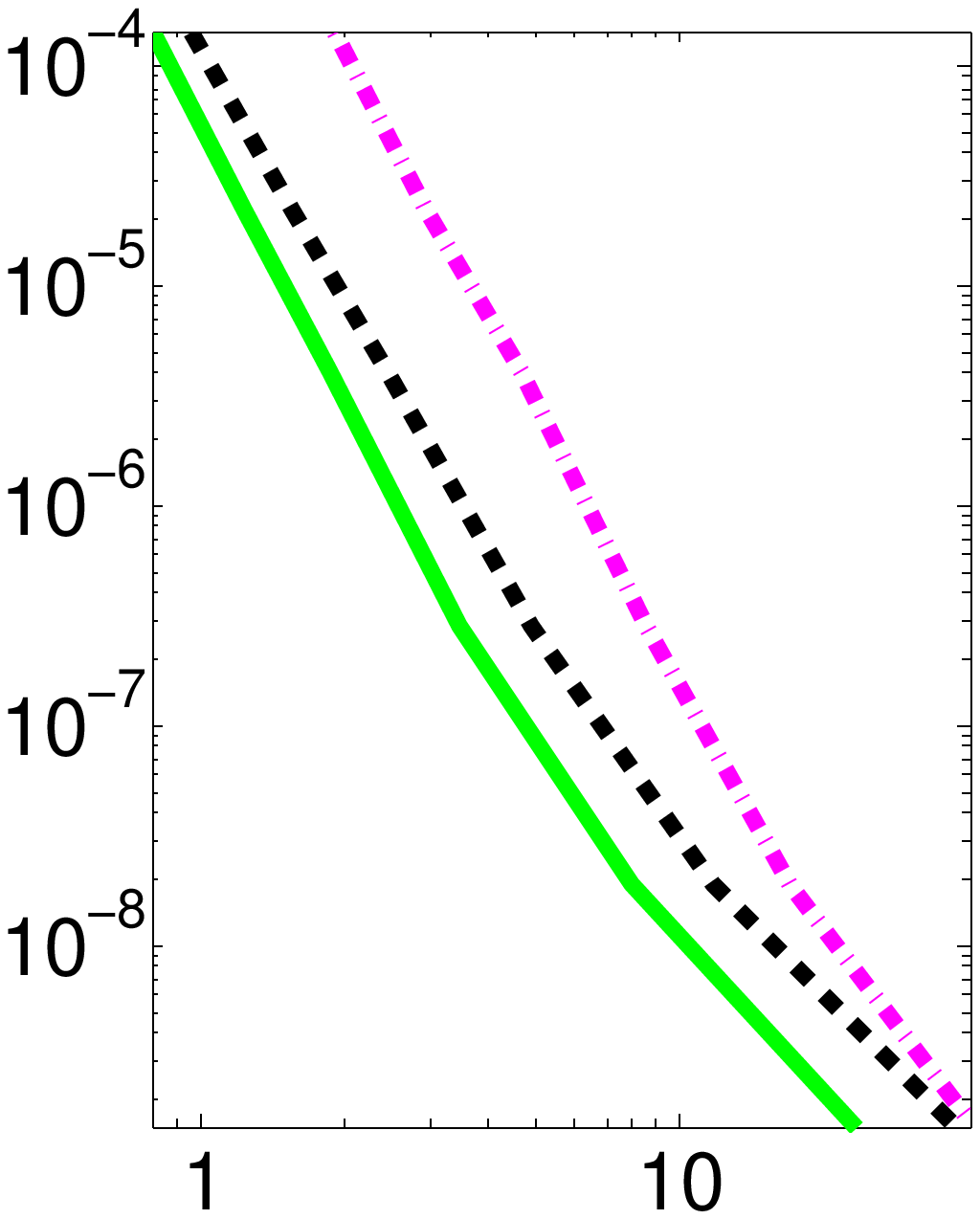}%
	    }}%
	    \;%
    	\vtop{\vskip0pt\hbox{%
	   		\includegraphics[width=\figsize]{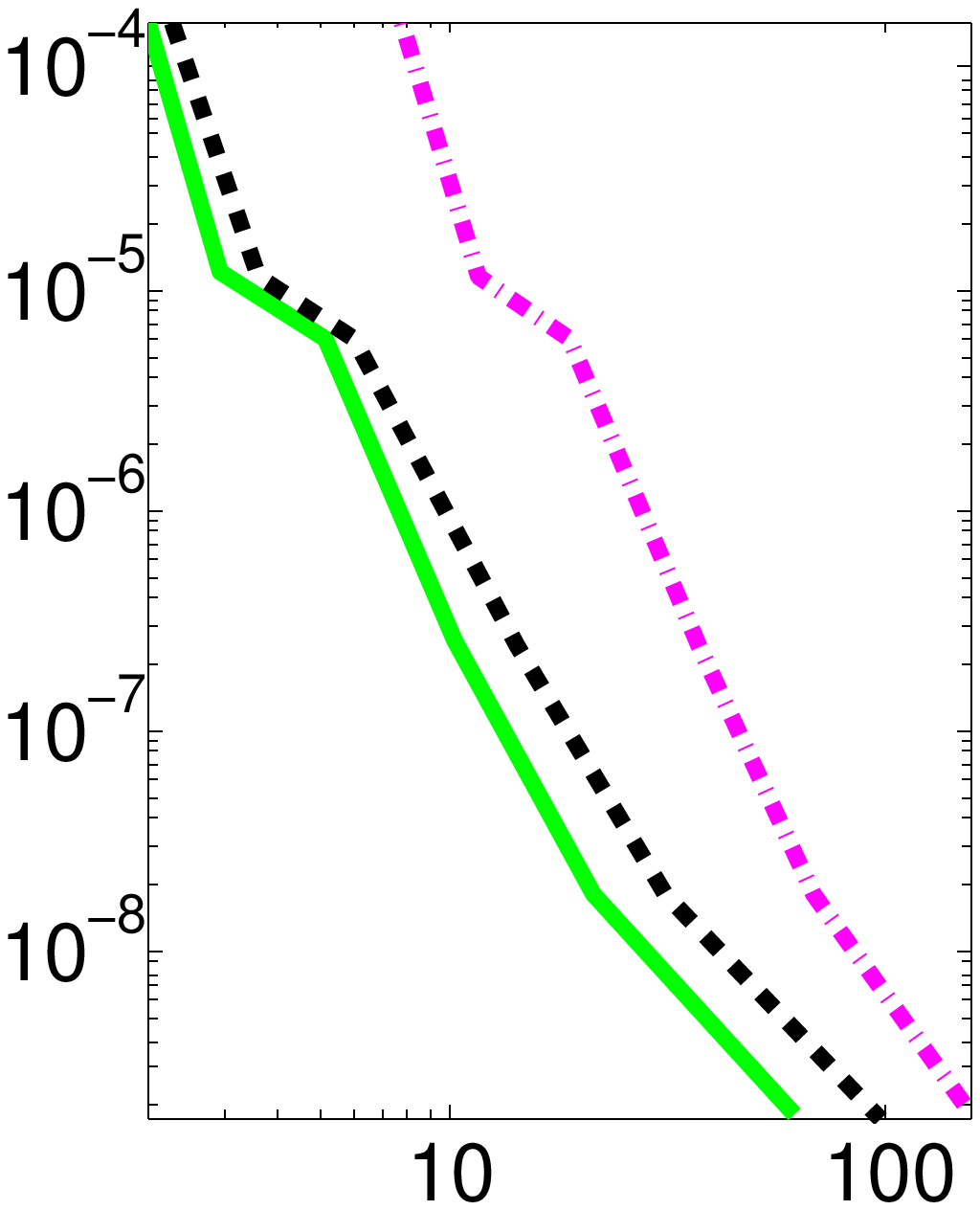}%
	    }}%
		\;%
    	\vtop{\vskip0pt\hbox{%
		    \includegraphics[width=\figsize]{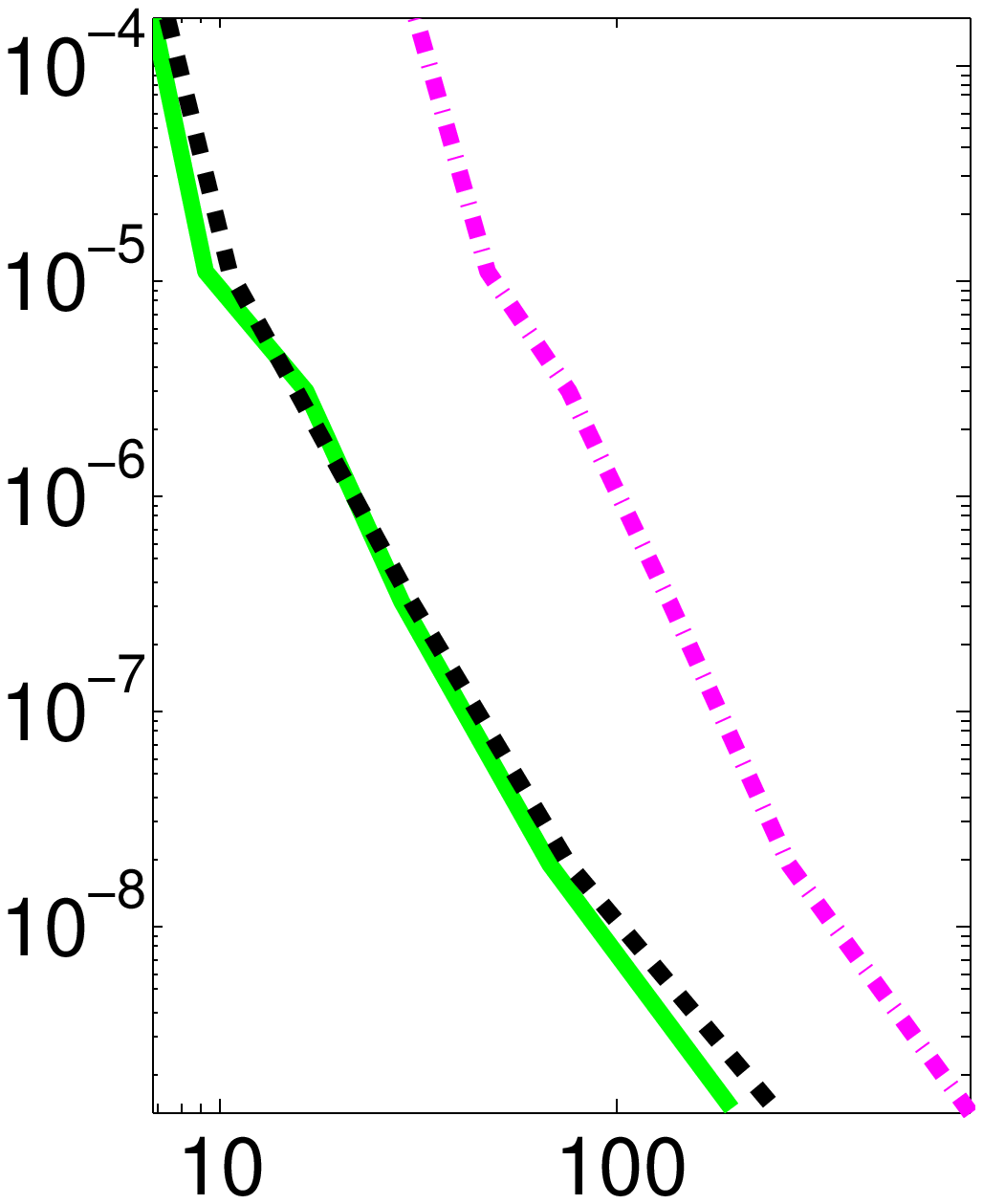}%
	    }}%
	}\\[0.25em]
    time (seconds) vs.\ error ($L^2$-Norm) \\[0.5em]
    \caption{Results for the Maxwell example. Integrators: 3 stage \EXPRB (solid line), \EXP4 (dashed line) and Runge-Kutta (dash dotted line), Ratios: \ratioA, \ratioB, \ratioC and \ratioD}
    \label{fig:MaxwellResults}
\end{figure}

The results are shown in Figure \ref{fig:MaxwellResults}. We give plots of spatial step size ($\Delta x$) vs.\ error (measured against the exact solution, not the reference solution), $\Delta x$ vs.\ chosen time step sizes ($h$) and a time vs.\ error plots.

In the first row of images we see, that the accuracy matches the requirements of the spatial discretization.

For the smallest grid ratio we see, that the exponential methods choose larger time steps. Still we get no speed gain, because we have to build an up to 36-dimensional Krylov subspace in each time step.

While both the Runge-Kutta and the exponential methods decrease their time step sizes due the bigger irregularity of the grid for shrinking grid ratios, the gap between the curves grows. This also corresponds to the time vs.\ error plots, where the exponential methods get relatively faster for more irregular grids.
In the following table we present the average quotient of the time step size of the \EXPRB integrator and the Runge-Kutta solver and the average speedup by grid ratio:
\centerline{
\begin{tabular}[t]{|c|c|c|c|c|}
\hline
ratio & \ratioA & \ratioB & \ratioC & \ratioD \\
\hline
avg. $h$ quotient (\EXPRB vs.\ Runge-Kutta) & \avgTimeStepQuotA & \avgTimeStepQuotB & \avgTimeStepQuotC & \avgTimeStepQuotD \\
\hline
avg. speedup (\EXPRB vs.\ Runge-Kutta) & \avgSpeedupA & \avgSpeedupB & \avgSpeedupC & \avgSpeedupD \\
\hline
\end{tabular}
}\\[0.5em]
We see that the exponential methods need to decrease their step size much less than the Runge-Kutta solver when decreasing the ratio.
such that they become relatively more effiecient on the more irregular grids. On the last grid \EXPRB is about five times faster than the conventional integrator.

The step size estimator was tested by manually increasing the time step size. In this case the solver became unstable, so the step size estimator actually detects the stability requirement of the method.

It should be noted that \EXPODE saves the solution in one long vector,
while the \MATLAB DG-codes use a custom format, where each field
component is saved as one rectangular matrix containing the degrees of
freedom for each finite element in its columns. For this reason our
\ode{} file has to switch between these two formats in each function
evaluation, which slows down the computation.

\subsection{Conclusions}

We presented a new advanced toolbox for exponential integrators, which implements several exponential integration schemes in recent research and uses modern techniques to approximate the arising matrix functions.
We designed the toolbox to solve large systems of differential equations, and we have shown its applicability to those systems here.
We also noticed rather weak performance for lower-dimensional systems, but the time cost was relatively small nonetheless.
This coincides with the package's philosophy: develop with small problems to use it with larger ones.

We allow substantial control over the integration process without the need of changing \EXPODE code directly.
See for instance the usage of the mass norm in Maxwell example.
This makes it easier to combine \EXPODE with other \MATLAB packages that perform well for other aspects of the equation.
For example, one may want to find a spatial discretization using
the DG-Codes from \cite{hesthaven2008nodal}.

It has been shown that our Krylov matrix function evaluator is able to
absorb spatial irregularity to some extent.
This lead to some ideas for the extension of the package by an improved Krylov method that can use the problem's mass matrix and mass scalar product to respect spatial structures.
We are currently implementing this idea, and the initial tests are quite promising.

The \EXPODE package is flexible, and we expect to implement further
exponential schemes in the future.

\section*{Acknowledgements}
	We thank Marlis Hochbruck for the initial idea for this package and for providing a lot of feedback.
	We also thank Julia Schweitzer for the initial \code{exprb}-code \EXPODE is built upon and Achim Sch\"adle for many great suggestions and a lot of valuable discussions.
	Patches provided by Abdullah Demirel to include order five exponential Rosenbrock-type schemes might be included in the next version.
	
	We thank the many testers, especially Anke Wortmann, Alexander Beckmann, Magdalena Weiss-Ribicky, Florian Kleen and everyone else who provided valuable feedback.

\bibliographystyle{plain}
\bibliography{expode}


\end{document}